\documentclass{article}

\usepackage{arxiv}

\usepackage[utf8]{inputenc} 
\usepackage[T1]{fontenc}    
\usepackage{hyperref}       
\usepackage{url}            
\usepackage{booktabs}       
\usepackage{amsfonts}       
\usepackage{nicefrac}       
\usepackage{microtype}      
\usepackage{lipsum}
\usepackage{graphicx}
\graphicspath{ {./figs/} }

\usepackage{amsmath}
\usepackage{bm}
\usepackage{xcolor}
\usepackage{subfig}

\def\xbold{\mathbf{x}}

\title{A coupled finite and boundary spectral element method for linear water--wave propagation problems

\thanks{\textit{\underline{Citation}}: 
\textbf{Antonio Cerrato, Luis Rodríguez-Tembleque, José A. González, M.H. Ferri Aliabadi. A coupled finite and boundary spectral element method for linear water-wave propagation problems. Applied Mathematical Modelling Volume 48, August 2017, Pages 1-20, DOI:10.1016/j.apm.2017.03.061.}} 
}

\author{
  Antonio Cerrato, Luis Rodríguez-Tembleque, José A. González  \\
  Escuela Ténica Superior de Ingeniería \\
  Universidad de Sevilla \\
  Camino de los Descubrimientos s/n, 41092 Sevilla, Spain\\
  \texttt{\{antoniocerrato, luisroteso japerez\}@us.es} \\
   \And
 M.H. Ferri Aliabadi \\
  Department of Aeronautics, Faculty of Engineering\\
   Imperial College of London\\
  South Kensington Campus, London SW7 2AZ, UK \\
  \texttt{m.h.aliabadi@imperial.ac.uk} \\
}

\begin{document}
\maketitle
\begin{abstract}
A coupled boundary spectral element method (BSEM) and spectral element method (SEM) formulation for the propagation of small-amplitude water waves over variable bathymetry regions is presented in this work. The wave model is based on the Mild-Slope Equation (MSE), which provides a good approach of the propagation of surface waves over irregular bottoms with slopes up to 1:3. In unbounded domains or infinite regions, the space can be divided into a central region or "inner region" and a surrounding infinite "outer region". The SEM allows us to model the inner region, where any variation of the bathymetry can be considered, while the outer region is modelled by the BSEM which, combined with the fundamental solution presented by Cerrato et al. (2016) \cite{Cerrato2016}, can include bathtymetries with straight and parallel contour lines. The solution approximation within the elements is constructed by a high order Legendre polynomials associated to Legendre-Gauss-Lobatto quadrature points, providing a spectral convergence for both methods. The proposed formulation has been validated in three different benchmark cases with different shapes of the bottom surface. The solutions present the typical $p$-convergence of spectral methods, with the combined advantages on the modelisation of the boundary element and finite element methods.

\end{abstract}


\section{Introduction}\label{sec1}


The pseudo-spectral methods have gained popularity over the last three decades since Patera \cite{Patera1984} presented his pioneer work on the mid 80s, combining the spectral approach with finite element formulation to solve the Navier-Stokes equation. Since then, the spectral element method has been used to solve many different problems in the fields of optics, electromagnetics, acoustics or water-wave propagation. In this area, much effort has been made for non-lineal equations like the Boussinesk-type equations \cite{Eskilsson2005,Eskilsson2006}. Nevertheless, the application of pseudo-spectral methodologies to the mild-slope equation has not been treated in the literature, despite of the fact that these type of models
represent the basic framework for the
simulation of surface wave propagation problems in variable water
depths.

The elliptical \textit{mild-slope} equation
(MSE) \cite{Berkhoff1972,Berkhoff1976} allows to consider simultaneously the effects of
diffraction, refraction, reflection and shoaling of linear water
surface waves, being formally only valid for slowly varying sea-bed
slopes, i.e. $\nabla h << kh$, being $h$ the water depth and $k$
the wave number \cite{Tsay1983,Booij1983}. That condition converts the MSE in a powerful tool to study large coastal areas with less computational effort compared to non-lineal models. In order to extend the applicability of the MSE, some modifications has been proposed for more complicated bathymetries, as in \cite{Porter1995} \cite{Massel1993}, or to include energy dissipation effects, such as wave breaking and bottom friction \cite{Maa2002}. The modified mild-slope equation (MMSE) presented in \cite{Chamberlain1995} and later improved in \cite{Porter1995}, retains the second order terms discarded by Berkhoff in the original formulation of the MSE. More extensions of the MSE can be found in \cite{Suh1997,Chandrasekera1997,Lee1998,Li1994,Hsu2000,Hsu2001}.

The MSE has been traditionally solved using the finite element
method (FEM) \cite{Berkhoff1976}, and the finite difference method (FDM) \cite{Panchang1991,Li1992}. However, classical finite difference and finite element methods present common deficiencies, as the pollution effect and difficulties to reproduce properly natural boundary conditions. 

The pollution effect is directly related to the numerical dispersion, which is an error due to the loss of the ellipticity of the Helmhotz's equation as the wave number $k$ increases. This error deteriorates the solution even if the coefficient $kh$, where $h$ is the element size, is kept small. Hence the number of nodes per wave length is not sufficient to determine the accuracy of the solution. This matter has been object of many works \cite{Bayliss1985,Ihlenburg1995-I,Ihlenburg1997-II,Deraemaeker1999}, where the pollution error is in general estimated as a rational function of $kh$ with the coefficients depending on the order of the interpolation function $p$. It is found that with higher $p$, the number of points per wave length required for a given accuracy presents a slower increasing rate.
In order to overcome this deficiency several formulations have been proposed, as the Galerkin Least-squares (GLS) methods \cite{Thompson1995}, spectral element methods \cite{Mehdizadeh2003} or the smoothed finite element methods \cite{Liu2006}. 

One of the most important problems of the FEM and FDM formulations is in the difficulty to reproduce unbounded domains properly. This deficiency has been studied by many authors. Bettess and
Zienkiewicz \cite{Bettess1978} and Lau and Ji \cite{Lau1989} used
infinite elements in the outer regions. Dirichlet to Neumann (DtN)
boundary conditions were proposed by Givoli et al.
\cite{Givoli1990,Keller1989,Givoli1991} as an analytical procedure
to reproduce exact non-reflecting boundary conditions in some
particular cases. This idea was followed by Bonet
\cite{Bonet2013} to derive the discrete non-local (DNL) boundary
condition. Other methods, such perfectly matched layers (PML) has been used in \cite{Giorgiani2013} and \cite{Modesto2015}. More rudimentary iterative methods have also been
proposed to define absorbing boundary conditions in \cite{Beltrami2001,Steward2001,Chen2002,Liu2008}.


Boundary element techniques prove to be very
accurate in wave-diffraction problems in unbounded
domains, presenting the additional benefit that the radiation
condition to infinity is automatically satisfied. In order to
improve the solution of the FEM schemes, Hauguel
\cite{Hauguel1978} and Shaw and Falby \cite{Shaw1978} proposed coupling FEM and BEM at the end of the seventies. Hamanaka \cite{Hamanaka1997} and later Isaacson and Qu \cite{Isaacson1990} introduced a boundary integral formulation to reproduce the wave field in harbors with partial reflecting boundaries and Lee et al. \cite{Lee2002,Lee2009} included the effect of incoming random waves. Boundary integral techniques have also been used for non-constant bathymetries by means of the dual reciprocity boundary element method (DRBEM). Some of these works also include wave run-ups, see Zhu \cite{Zhu1993}. Later, this technique was extended to model internal regions with variable depth surrounded by exterior
regions with constant bathymetry \cite{Liu2003,Zhu2000,Zhu2009,Hsiao2009}. More recently, Naserizabeh et al. \cite{Naserizadeh2011} proposed a coupled BEM-FDM formulation to solve the MSE in unbounded domains. Later, Cerrato et al. \cite{Cerrato2016} derived a complete kernel of a fundamental solution for variable water depths, based on the Green's function presented in \cite{Belibassakis2000}, and introduced it in a standard BEM formulation.


In the context of SEM, we have the contributions of Mehdizadeh and Paralchivoiu (2003) \cite{Mehdizadeh2003} or, more recently, He et al. (2016) \cite{He2016}. These works presented a SEM formulation based on the Legendre polynomials for the two dimensional Helmholtz's equation using special techniques to reproduce unbounded domains. To model the open boundaries, Mehdizadeh and Peralchivoiu use an absorbing layer, or PML, surrounded the inner region while He et al. used the Dirichlet to Neumann (DtN) boundary conditions. This last work is also concerned with the modelling of waves propagating over layered media. An specific work about water wave modelling with BEM combined with spectral approximation based on Chebyshev polynomials was presented recently by Kumar et al. \cite{Kumar2015}, but their formulation is limited to constant bathymetries.
The combination of spectral elements and BEM formulations has been also recently applied to other fields, such as elastodynamic promblems \cite{Zou2016}.
A final interesting work is \cite{Vos2010}, where the authors study the $hp$-convergence of the spectral element method for the two dimensional Helmholtz's equation.

In order to deal with the issues mentioned above, this work presents a spectral boundary element-finite element formulation for linear water-wave propagation problems. The SEM is used to model inner regions, allowing to model bathymetries with arbitrary shape, while the BSEM is used as a natural boundary condition which fulfils automatically the Sommerfeld radiation condition. An important feature of the BSEM formulation is that it uses the fundamental solution derived by Cerrato et al. \cite{Cerrato2016} to consider a variable bathymetry in the outer region. The approach of the solution is made by considering a pseudospectral approximation inside the elements, using nodal basis functions with nodes located at the LGL points. 

The paper is organized as follows. First, Section \ref{sec2} presents the MSE for the linear water-wave propagation problem. In Section \ref{sec3}, the boundary spectral element formulation is explained in detail. In Section \ref{sec4}, special attention is paid to the finite spectral element formulation. The coupled boundary spectral element--finite spectral element formulation is presented in Section \ref{sec5}. Section \ref{sec6} is dedicated to the validation of
the proposed formulation through the solution of different wave
propagation problems for variable water depth, where mathematical and numerical convergence studies are presented. Finally, the paper concludes with the summary and conclusions.

\section{The mild-slope equation}\label{sec2}
To formulate the MSE, according to
\cite{Berkhoff1972,Berkhoff1976}, we consider a Cartesian
coordinate system with the \((x,y)\)-plane located on the quiescent
water surface and the $z$ direction pointing upwards. The still
water depth is given by $h_w(x,y)$ and $\nabla=(\partial_{x},
\partial_{y})$ is used to represent the gradient operator.
Under the assumption of potential flow and integrating the
velocity potential in the vertical direction with appropriated
boundary conditions, the velocity potential of the water surface
can be assumed to be of the form:
\begin{equation}
\Phi(x,y,t)=\phi(x,y) e^{-i\omega t},
\end{equation}
where $i$ is the imaginary unit and $t$ is the time variable. This
potential has to satisfy the homogeneous MSE, that may be written
as:
\begin{equation}\label{eqn-MSE}
\nabla \cdot (cc_{g}\nabla \phi) + \omega \dfrac{c_{g}}{c} \phi =
0.
\end{equation}
where $c$ is the wave velocity and $c_{g}$ is
the group velocity. The water depth function $h_w(x,y)$, the wave
number $k$ and the angular frequency $\omega$ of the waves are
related by the dispersion equation:
$\label{eqn:dispersionEquation} \omega^{2}=gk \tanh(kh_w)$, where
$g$ is the gravitational acceleration ($g=9.81m/s^{2}$). So, for a
fixed frequency and variable bathymetry, the wave number $k(x,y)$
is a function of the local water depth. In the linear wave theory, the wave height ($H$), is linearly related to the velocity potential on the water surface by the following relation:
\begin{equation}
H=2\omega|\phi|/g,
\end{equation}
an important design variable for practical problems.

The MSE can be simplified introducing the following change of
variable due to Bergmann \cite{Bergmann1946}:
$\phi={\hat{\phi}}/{\sqrt{cc_{g}}}$. This relation transforms
\eqref{eqn-MSE} into a Helmholtz equation:
\begin{equation}
\label{eqn:Helmholtz} \nabla^{2} \hat{\phi} + \hat{k}^{2}
\hat{\phi} = 0,
\end{equation}
with a modified wave number $\hat{k}(x,y)$ given by:
\begin{equation}
\label{eqn:k con gorro} 
\hat{k}^{2}(x,y)= k^{2} -
\frac{\nabla^{2}\sqrt{cc_{g}}}{\sqrt{cc_{g}}},
\end{equation}
that is a known function of the wave characteristics and the local
water depth. In this work, we discard the second term and the modified wave number is approached as $\hat{k}=k(x,y)$.

Note that a similar change of variable can be done for treating the same problem under the framework of the MMSE. Simply by modifying accordingly the expression of the wave number \eqref{eqn:k con gorro}, to account for additional effects associated with higher-order contributions of bottom slope and curvature, we obtain a MMSE model that extends the applicability of the mild-slope equation.

\section{Boundary spectral element method}\label{sec3}
The Boundary Integral Equation (BIE) is a powerful tool to construct accurate and efficient formulations for linear water wave problems. The BIE presents some important advantages as a the fullfilment of the Sommerfeld radiation condition at infinity, or the strong reduction on the number of degrees of freedom required to construct the discrete model of the domain. The main disadvantage of the application of the BIE to water wave propagation problems is the difficulty to include variable bathymetries. The prefered technique to deal with this difficulty is the DRBEM. However it requires to compute domain integrals and the modelisation of infinite domains continues being only applicable to constant-depth areas. Nevertheless, in order to model infinity domains with non-homogeneus water-depth, we can use the standard BIE but equipped with an appropriated fundamental solution \cite{Cerrato2016}. 

In the BSEM, the BIE is discretized by spectral elements. As in other spectral element techniques, the nodal basis functions based on a family of orthonormal polynomials provide to the BIE an exponencial convergence and highly accurate description of the boundary shape.

\subsection{Boundary integral equation}\label{sec31}

We start writing the classical BIE of the BEM in order to fix the notation. Considering a boundary $\Gamma$ of the region of interest, which is governed by equation \eqref{eqn:Helmholtz}, the BIE can be expressed as:
\begin{equation}
C(\mathbf{x'}) \hat{\phi}(\mathbf{x'})  +
\int_{\Gamma} \nabla \psi(\mathbf{x},\mathbf{x'};\hat{k}) \cdot \mathbf{n}~\hat{\phi}(\mathbf{x}) \:d\Gamma
=
\int_{\Gamma}  \psi(\mathbf{x},\mathbf{x'};\hat{k})~\hat{q}(\mathbf{x}) \:d\Gamma + \hat{\phi}_{in},
\label{eqn:BEMBoundaryIntegralEquation}
\end{equation}
where $\mathbf{x'}$ is the collocation point, $\hat{q}$ is the normal flux through the boundary considering an outward normal and $\psi(\mathbf{x},\mathbf{x'};\hat{k})$ is the Green's function of the Helmholtz problem. The free term is, in general, given by $C(\mathbf{x'})={\theta(\mathbf{x'})}/{2\pi}$, being $\theta(\mathbf{x'})$ the internal angle of the boundary at the collocation point. In order to consider also scattering problems, the incident potential $\hat{\phi}_{in}$ is included in the above expression \cite{Wu2000}.

\subsection{Spectral boundary elements}\label{sec32}

\begin{figure}
\fontsize{8}{8}\selectfont 
  \centering  
  \def\svgwidth{0.8\columnwidth}
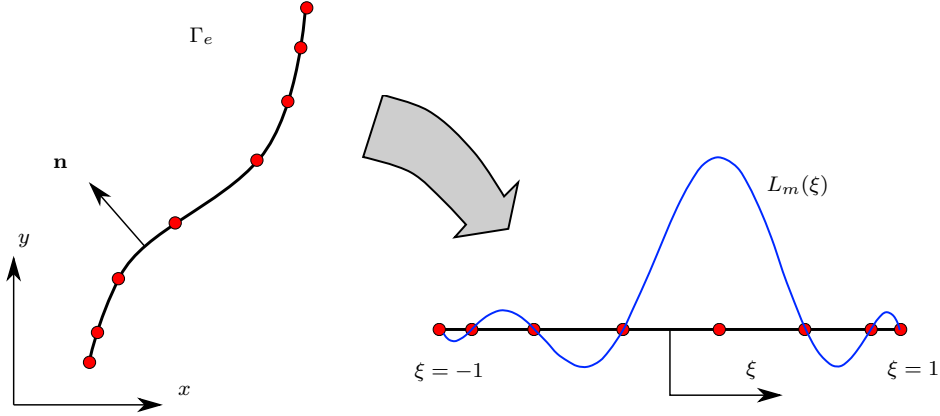
\normalsize
\caption{Example of a spectral boundary element. An element of order $p=7$ is represented in the global $x-y$-axis (left) and in the local coordinate $\xi$ (rigth). Also the nodal function corresponding to the fifth node is shown} 
\label{fig:Spectral boundary element}
\end{figure}
%

In order to evaluate the integrals of equation \eqref{eqn:BEMBoundaryIntegralEquation} the boundary is divided into $n_e$ elements, being the discrete boundary $\Gamma = \bigcup_{i=1}^{n_e} \Gamma_{e}$. In the BSEM formulation presented here, the geometry and the variables on the boundary are approximated using nodal basis functions based on the Lagrange polynomials with nodes located on the Legendre-Gauss-Lobatto (LGL) quadrature points. Inside each spectral element, the modified velocity potential and the normal flux through the boundary is approximated by: 
\begin{eqnarray}\label{eqn:BEMdiscretizacion1}
\hat{\phi}(\xi)=\sum_{m=0}^{p}
L_{m}(\xi) \hat{\phi}_m^{(e)},
& &
\hat{q}(\xi)=\sum_{m=0}^{p}
L_{m}(\xi) \hat{q}_m^{(e)},
\end{eqnarray}
where $\hat{\phi}_{m}^{(e)}$ and $\hat{q}_{m}^{(e)}$ are the nodal
values of the element $e$ and $L_{m}(\xi)$ the nodal basis functions of degree $p$ on the interval $\left[-1,1\right]$. This approximation functions can be constructed by means of Lagrange polynomials with roots at the Gauss-Lobatto points:
\begin{equation}\label{eqn:Lagrange_interpolant_1D}
L_{m}(\xi)= \prod\limits_{l=0 (l\neq m)}^{p}\frac{\xi - \xi_l}{\xi_m - \xi_l}, \quad m=0,1,...,p,
\end{equation}
being $\xi_m$ the local coordinate of the node under consideration and $\xi_l$ the coordinate of the rest of the nodes inside the element. In the above expressions the index $0$ and $p$ are used to denote the nodes located at the extremes of the interval ($\xi=-1$ and $\xi=1$). A representation of this spectral boundary element is given in Figure \ref{fig:Spectral boundary element}. It is possible to use the same high order polynomials to approximate the geometry, specially  when the real geometry is curved and complex.

\subsection{Discrete boundary integral equation}\label{sec33}
Introducing in the BIE given by \eqref{eqn:BEMBoundaryIntegralEquation} the approximations of the variables \eqref{eqn:BEMdiscretizacion1}, the discrete form of the BIE can be written as
\begin{multline}
\label{eqn_BEMdiscretized BIE}
C_i\hat{\phi}_i +  \sum_{e=1}^{n_{e}} \int_{-1}^{1} \nabla \psi(\mathbf{x}(\xi),\mathbf{x}_i;\hat{k}) \cdot \mathbf{n} ~L_{m}(\xi) J^{(e)}(\xi)  \: d\xi \: \hat{\phi}_{m}^{(e)}
=
\\
\sum_{e=1}^{n_{e}} \int_{-1}^{1}  \psi(\mathbf{x}(\xi),\mathbf{x}_i;\hat{k}) ~ L_{m}(\xi) J^{(e)}(\xi)  \: d\xi \:  \hat{q}_{m}^{(e)}
+
\hat{\phi}_{in\:i}
\end{multline}
where the sub-index $i$ refers to the collocation point and the super-index $e$ to the elements. The jacobian of the transformation from a global coordinate to the local coordinate system is denoted by $J^{(e)}(\xi)$. Assembling all the element contributions included in expression \eqref{eqn_BEMdiscretized BIE}, the discrete form of the BIE can be finally written in a matrix form as:
\begin{equation}\label{eqn:BSEM_system}
\mathbf{\hat{H}} \boldsymbol{\hat{\phi}} =
\mathbf{\hat{G}}\mathbf{\hat{q}} + \boldsymbol{\hat{\phi}}_{in},
\end{equation}
where the vectors $\boldsymbol{\hat{\phi}}$ and $\mathbf{\hat{q}}$ contain the unknown nodal values of the modified potential and flux respectively, and $\boldsymbol{\hat{\phi}}_{in}$ is the load vector with the values of the incident wave. The matrices $\mathbf{\hat{H}}$ and $\mathbf{\hat{G}}$ are composed of the following terms:
\begin{align}
\label{eqn:BSEM_Hij}
\hat{H}_{ij} &= C_i \delta_{ij} + 
  \sum_{e=1}^{E^{(j)}} \int_{-1}^{1} \nabla \psi(\mathbf{x}(\xi),\mathbf{x}_i;\hat{k}) \cdot \mathbf{n} ~ L_{j}(\xi) J^{(e)}(\xi) \:d\xi,
\\
\label{eqn:BSEM_Gij}
\hat{G}_{ij} &= \sum_{e=1}^{E^{(j)}} \int_{-1}^{1} \psi(\mathbf{x}(\xi),\mathbf{x}_i;\hat{k}) ~ L_{j}(\xi) J^{(e)}(\xi) \:d\xi,
\end{align}
where $E^{(j)}$ represents the group of elements sharing node $j$.

\subsection{Numerical integration}

The integrals of the $\hat{H}_{ij}$ and $\hat{G}_{ij}$ terms are evaluated numerically. One of the principal advantages of using a Lagrangian approach combined with a LGL distribution of the nodes within the element, is that the numerical integration can be carried out by means of its associated LGL quadrature rule, which is exact for polynomials up to $2p-1$. This is specially accurate for those elements that are not located close to the collocation point, over a region with smooth variation of the wave number and without a complex geometry. In this case, the values of the Lagrange interpolants at the quadrature points are given by the relation 
\begin{equation}
\label{eqn:Lxidelta}
L_m(\xi_k)=\delta_{mk},
\end{equation}
being $\delta$ the Kronecker delta.
Using this property, the integral of the product of a function $f(\xi)$ with the nodal basis functions $L_m(\xi)$ reduced to a single evaluation of the integrand at the quadrature point as follows:
\begin{equation}
\label{eqn:BSEM_numerical_integration}
\int_{-1}^{1} f(\xi)L_m(\xi)\:d\xi 
 = f(\xi_m) w_m,
\end{equation} 
where the weight is given by the expression:
\begin{equation}
\label{eqn:LGLweights}
w_m=\dfrac{2}{p(p+1)} \dfrac{1}{P_p^2(\xi_m)},
\end{equation}
and where $P_p(x)$ is the Lagrange polynomial of order $p$. This way, the values of the integrals inside \eqref{eqn:BSEM_Hij} and \eqref{eqn:BSEM_Gij} can be computed fast and efficiently. 

In other cases, like for example integration on elements close to the collocation point, over regions where the geometry contains a sharp corner, or the wave number changes with a complex function, a standard Gaussian quadrature is preferred. 

\subsection{Regularization of the singular integrals}\label{sec34}

\begin{figure}
\fontsize{8}{8}\selectfont
\centering
\def\svgwidth{0.8\columnwidth}
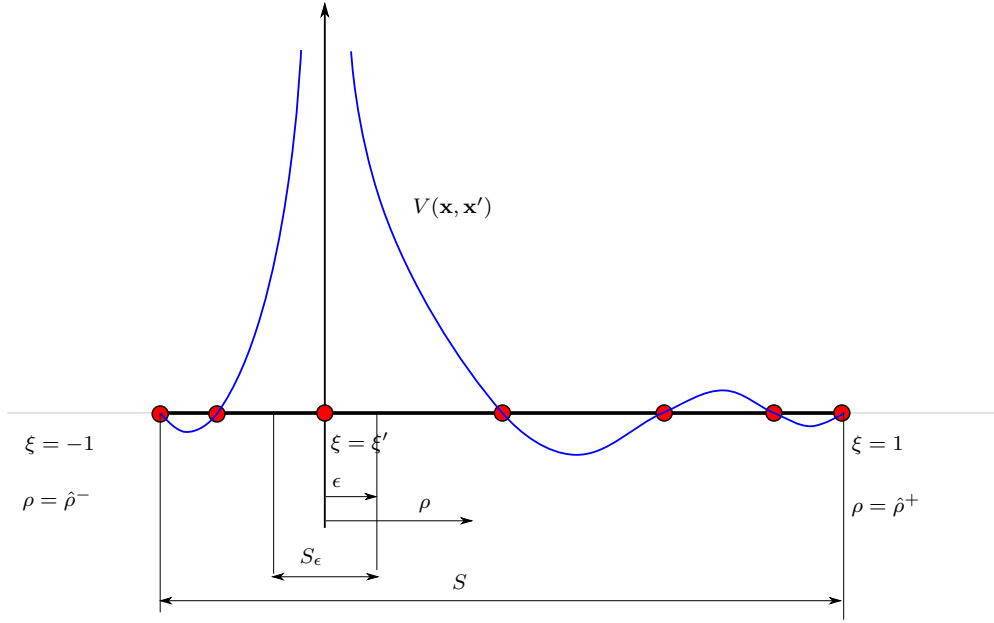
\normalsize
\caption{Singular integral around $S_{\epsilon}$ when the collocation point is located inside the element. The distance to the singularity $\rho=|\xi-\xi'|$ is defined to extract the singular kernel}
\label{fig:Sing integral}
\end{figure}
%

During the collocation process, when the element under integration contains the source load, the integrand of \eqref{eqn:BSEM_Gij} becomes weakly singular. This problem can be treated by a regularisation technique; subtracting the singularity and integrating it separately \cite{Aliabadi1985,Aliabadi1989,Guiggiani1992}. 

Let us define the $S$ and the $S_\epsilon$ regions inside the element as shown in Figure \ref{fig:Sing integral}, where $S_\epsilon$ surrounds the singularity and $S$ covers the entire element. Taking the limit when $\epsilon$ tends to zero, the singular integrals in \eqref{eqn:BSEM_Gij} can be calculated by means of their Cauchy Principal Value:
\begin{equation}
\label{eqn: Apend Sing Int I}
I=\lim_{\epsilon\longrightarrow 0} 
\left\lbrace
 \int_{S-S_{\epsilon}} V(\xbold,\xbold')- V_{1}(\xbold,\xbold') \: d\xi \right.
 \: + \: 
\left. \int_{S-S_{\epsilon}}  V_{1}(\xbold,\xbold') \: d\xi 
\right\rbrace
\end{equation}
where the integrands are 
$V(\xbold,\xbold') = \psi(\xbold(\xi),\xbold(\xi')) L_i(\xbold(\xi)) J(\xbold(\xi))$ 
and 
$V_{1}(\xbold,\xbold')$ is the first term of the Taylor series expansion of $V(\xbold,\xbold')$ at $\xi=\xi'$.

To subtract the singularity, we use the fundamental solution of the 2D-Helmholtz equation for constant wave number in terms of its series expansion: 
\begin{equation}
\label{eqn: Apend Sing Int FS approach}
\psi(\xbold,\xbold')= \wp(k\sigma) +\dfrac{i}{4} +
O((k\sigma)^2 log(k\sigma)),
\end{equation}
where 
\begin{equation}
\label{eqn: Apend Sing}
\wp(k\sigma)=-\dfrac{1}{2\pi}(log(k\sigma/2)+\gamma)
\end{equation}
and $\sigma=|\xbold-\xbold'|$ is the distance to the load. Then, in order to express the weakly singular kernel of \eqref{eqn: Apend Sing Int I} when $\sigma$ goes to zero, the function containing the singularity \eqref{eqn: Apend Sing} is written in terms of the local variable of the element $\rho=\vert \xi-\xi' \vert$, as shown in Figure \ref{fig:Sing integral}. Then, the distance from $\xbold'$ in any direction $x_i$ can be approximated using the first order term of its Taylor series expansion: 
\begin{equation}
\label{eqn: Apend FOA xi}
x_i-x'_{i}= \left. \dfrac{dx_i}{d\xi} \right|_{\xi=\xi'} \rho + O(\rho^2).
\end{equation}
Therefore, the distance between $\xbold'$ and any point $\xbold$ can be finally expressed as
\begin{equation}
\label{eqn: Apend FOA r}
\sigma=\rho A + O(\rho^2),
\end{equation}
where $A^2= \sum_{i=1}^{2} \left( \left. x_{i,\xi} \right|_{\xi=\xi'}\right) ^2$ contains the Jacobian of the transformation. 

Substituting \eqref{eqn: Apend FOA r} into \eqref{eqn: Apend Sing} we obtain the kernel of the singularity in terms of the local distance $\rho$:
\begin{equation}
\wp(k\sigma)=\wp_o+\wp_1 \log{\rho} + O(\rho),
\end{equation}
where $\wp_o=-\dfrac{1}{2\pi} (\gamma + \log{\dfrac{kA}{2}})$ and $\wp_1 = -\dfrac{1}{2\pi}$.

The same procedure can be applied to all the terms appearing in the integrand, like:
\begin{equation}
L(\rho)=L_o + O(\rho)
\end{equation} 
and
\begin{equation}
J(\rho)=J_o + O(\rho),
\end{equation} 
being $L_o=L_i(\xi')$ and $J_o=J(\xi')$. Collecting all the terms, we can write:
\begin{equation}
V(\xbold,\xbold')
=
T_o + T_1 \log{\rho} + O(\rho)
\end{equation} 
where $T_o=\wp_o L_oJ_o$ and $T_1=\wp_1 L_oJ_o$. Therefore, the first term of the series expansion of the integrand $V(\xbold,\xbold')$ can be written as:
\begin{equation}
\label{eqn: Apend Integrand approach}
V_1(\xbold,\xbold')=T_o + T_1 \log{\rho}.
\end{equation}
Introducing this in equation \eqref{eqn: Apend Sing Int I} we regularize integral in the domain $S_{\epsilon}$. 

To compute the first regular integral of \eqref{eqn: Apend Sing Int I} we can use a standard Gaussian quadrature. The second integral, which contains the singularity, can be computed analytically dividing the domain $S$ into two regions:
\begin{equation}
\int_{S} V_1(\xbold,\xbold') \:d\xi 
= 
\int_{\hat{\rho}^{-}}^{0}  T_o + T_1 \log{\rho} \: d\rho
+
\int_{0}^{\hat{\rho}^{+}}  T_o + T_1 \log{\rho} \: d\rho,
\end{equation}
and integrating,
\begin{equation}
\int_{S} V_1(\xbold,\xbold') \:d\xi 
= 
2T_o+T_1 \left( \hat{\rho}^{-} \log{\hat{\rho}^{-}} + \hat{\rho}^{+} \log{\hat{\rho}^{+}} -2 \right),
\end{equation}
being $\hat{\rho}^{-}=|-1-\xi'|$ and $\hat{\rho}^{+}=|1-\xi'|$.

\section{Spectral element method}\label{sec4}

The SEM is a high-order finite element
technique that combines the geometrical flexibility of finite
elements with the high accuracy of spectral methods. This method
was pioneered in the mid 80's by Patera \cite{Patera1984} for the Navier-Stokes equation 
and exhibits several favorable computational
properties, such as the use of naturally diagonal
mass matrices that facilitate the iterarive solution process. Moreover, in the context of Helmholtz's equation, the SEM reduces the pollution effect and improves the efficiency of the classical Galerkin finite element approximation, i.e., it requires a lower number of nodes per wavelength compared to FEM  in order attain the same accuracy \cite{Mehdizadeh2003}.

\subsection{Weak form of the MSE}

The weighted residual form of the Helmholtz's equation  \eqref{eqn:Helmholtz} for a domain $\Omega$ bounded by $\Gamma$ after integrating by parts is expressed as:
\begin{equation}
\label{eqn:Weighted Residual}
\int_{\Omega} \nabla\eta \nabla\hat{\phi}\:d\Omega -
\int_{\Gamma}  \eta\hat{q} \:d\Gamma - \int_{\Omega} \hat{k}^2
\eta\hat{\phi}\:d\Omega
=
0, 
\end{equation}
where $\eta$ is the weight function. This variational formulation is the starting point of the SEM.

\subsection{Two dimensional spectral elements}
\begin{figure}[t]
\fontsize{8}{8}\selectfont
\centering
\def\svgwidth{10cm}
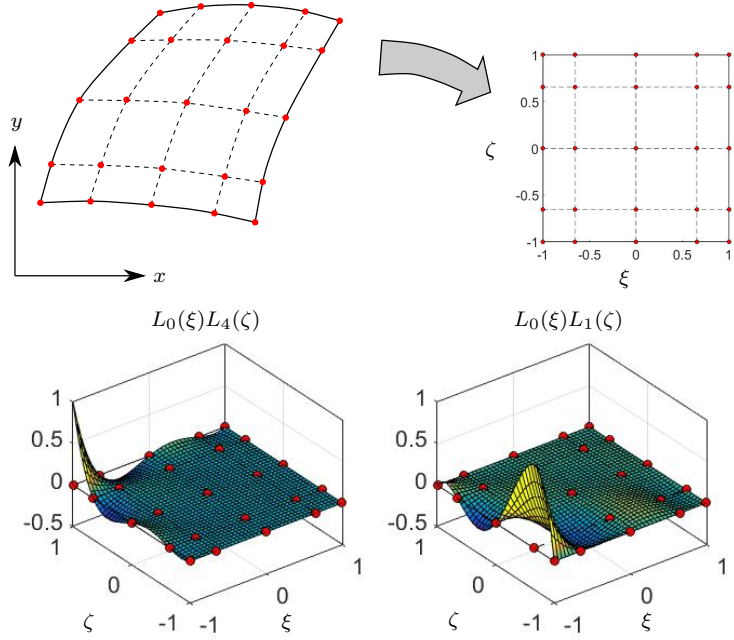
\normalsize
\caption{A finite spectral element of order $p=4$ is represented in the physical ($x-y$)-coordinate system and in the normalized $(\xi,\zeta)\in[-1,1]^2$ reference coordinate system (top). Two representative nodal basis functions of the element are also showed (bottom)} \label{fig:Spectral finite element}
\end{figure}

A discretization of the domain $\Omega$ is carried out by using $n_e$
spectral elements, fulfilling
\begin{equation}
\label{eqn:domain decomposition}
\Omega = \bigcup_{e=1}^{n_e} \Omega_{e} 
\quad
\text{and}
\quad
\bigcap_{e=1}^{n_e} \Omega_{e}=\emptyset,
\end{equation}
where the modified velocity potential is approached inside the elements by means of Lagrangian interpolants through the LGL grid points in the form:
\begin{equation}
\label{eqn:FEMdiscretizacion1} \hat{\phi}(\xi,\zeta)=\sum_{m=0}^{p}
\sum_{n=0}^{p} L_{m}(\xi)L_{n}(\zeta)  \hat{\phi}(\xi_m,\zeta_n).
\end{equation}
where $L_m$ is the one dimensional shape function associated to node $m$ as in \eqref{eqn:Lagrange_interpolant_1D}, $p$ is the polynomial degree and $(\xi,\zeta)\in[-1,1]^2$ are the intrinsic coordinates of the reference element. Figure \ref{fig:Spectral finite element} shows an example of a two-dimensional spectral element and its transformation to the normalized coordinates $(\xi,\zeta)$. Also two the nodal basis functions corresponding to $\mathcal{L}_{21}(\xi,\eta)=L_0(\xi)L_4(\zeta)$ and $\mathcal{L}_{6}(\xi,\eta)=L_0(\xi)L_1(\zeta)$ are represented. 

Equation \eqref{eqn:FEMdiscretizacion1} can be written in a more compact form as
\begin{equation}
\label{eqn:FEMdiscretizacion2}
\hat{\phi}(\boldsymbol{\xi})=\sum_{j=1}^{N_c} \mathcal{L}_{j}(\boldsymbol{\xi})
\hat{\phi}_j,
\end{equation}
where $\boldsymbol{\xi}=(\xi,\zeta)$, $N_c$ is the total number of element nodes ($(p+1)^2$ for quadrilateral elements),
$\mathcal{L}_{j}(\boldsymbol{\xi})$ the two dimensional shape functions and $\hat{\phi}_j$ the nodal values of the modified velocity potential. 

%

\subsection{Discrete finite element equations}

Applying the Spectral-Galerkin method to the weighted residual form \eqref{eqn:Weighted Residual} and considering the domain discretization \eqref{eqn:domain decomposition}, the following  system of equation is obtained:
\begin{eqnarray}\label{eqn:SEM_system}
\left[
\begin{array}{ccc}
\hat{\mathbf{A}} & -\hat{\mathbf{C}}
\end{array}
\right]
\left\{
\begin{array}{c}
    \hat{\boldsymbol{\phi}}
    \\
    \hat{\mathbf{q}}
\end{array}
\right\} 
\!\!=  
\left\{\begin{array}{c}
    \mathbf{0} 
    \end{array}
    \right\},
\end{eqnarray}
where vectors $\hat{\boldsymbol{\phi}}$ and
$\hat{\mathbf{q}}$ group all the nodal values of the modified velocity potentials and boundary fluxes respectively. The matrix $\hat{\mathbf{A}}$ and $\hat{\mathbf{C}}$ are obtained by assembling the elemental matrices $\hat{\mathbf{A}}^e$ and $\hat{\mathbf{C}}^e$ as in the standard finite element formulation. The domain integrals of \eqref{eqn:Weighted Residual} over each element $e$ provide the coefficients of the $\hat{\mathbf{A}}^e$ matrix:
\begin{equation}
\hat{A}_{ij}^e = K_{ij}^e - M_{ij}^e,
\end{equation} 
where
\begin{equation}
\label{eqn:Stiffness matrix}
K_{ij}^e = 
\int_{\Omega_{e}} \nabla \mathcal{L}_{i}\nabla \mathcal{L}_{j} \:d\Omega
\end{equation}
and
\begin{equation}
\label{eqn:Mass matrix}
M_{ij}^e = 
\int_{\Omega_{e}}  \hat{k}^2(x,y)\mathcal{L}_{i}\mathcal{L}_{j}  \:d\Omega.
\end{equation}
The coefficients of the $\hat{\mathbf{C}}^e$ matrices correspond to the boundary integrals of \eqref{eqn:Weighted Residual}, being:
\begin{equation}
\label{eqn:Boundary matrix}
\hat{C}_{ij}^{e}  = \int_{\Gamma_{e}} L_{i} L_{j} \:d\Gamma. 
\end{equation}

\subsection{Numerical integration}

To compute the integrals appearing in the SEM, the LGL quadrature is selected due to the convenient evaluation of the Lagrange interpolants at the quadrature points, as it was shown in \eqref{eqn:Lxidelta}. In the two dimensional form, the quadrature is written as
\begin{equation}
\label{eqn:SEM_numerical_integration}
\int_{-1}^{1} \int_{-1}^{1} f(\boldsymbol{\xi})\:d\xi \:d\zeta
 = \sum_{j=1}^{N_c}  f(\boldsymbol{\xi}_j) \bar{w}_j,
\end{equation} 
for a general integrand $f(x)$, where $\boldsymbol{\xi}_j$ are the local coordinates of each quadrature point $j$ and $\bar{w}_j$ are their corresponding weights. 

Hence, the computation of the "stiffness" matrix \eqref{eqn:Stiffness matrix} by the LGL quadrature of order $p$ leads to a diagonal banded matrix. Moreover the integrand is composed by polynomials of order $2p-2$ on each direction, being the LGL quadrature exact for polynomials of order up to $2p-1$. 

One of the main advantages of the LGL quadrature of order $p$ is that the domain integrals \eqref{eqn:Mass matrix} and the boundary integrals \eqref{eqn:Boundary matrix} produces diagonal matrices. The coefficients are easily obtained merely by evaluating the jacobian at the quadrature points and multiplying it by their corresponding weights. Therefore, the "mass" matrices are obtained by
\begin{equation}
\int_{-1}^{1} \int_{-1}^{1} 
\hat{k}^2(\boldsymbol{\xi})
\mathcal{L}_i(\boldsymbol{\xi}) \mathcal{L}_j(\boldsymbol{\xi})
J(\boldsymbol{\xi}) \:d\xi \:d\zeta
\simeq
\hat{k}^2(\boldsymbol{\xi}_i) J(\boldsymbol{\xi}_i) \bar{w}_i \delta_{ij},
\end{equation}
and the boundary integrals
\begin{equation}
\int_{-1}^{1} L_i(\xi) L_j(\xi) J(\xi) \:d\xi 
\simeq 
J(\xi_i) w_i \delta_{ij},
\end{equation}
which are expressions easy and efficient to calculate. However these integrals are computed with some error due to the limited accuracy of the quadrature but, the spectral convergence of the SEM is maintained even for non-homogeneous medium or complex geometries \cite{Maday1990}. Nevertheless a standard quadrature is preferred for some cases, as for example when the wave number function $\hat{k}(x,y)$ presents a non-continuous derivative.

\section{BSEM and SEM coupling}\label{sec5}
\begin{figure}[t]
\centering
\def\svgwidth{0.8\columnwidth}
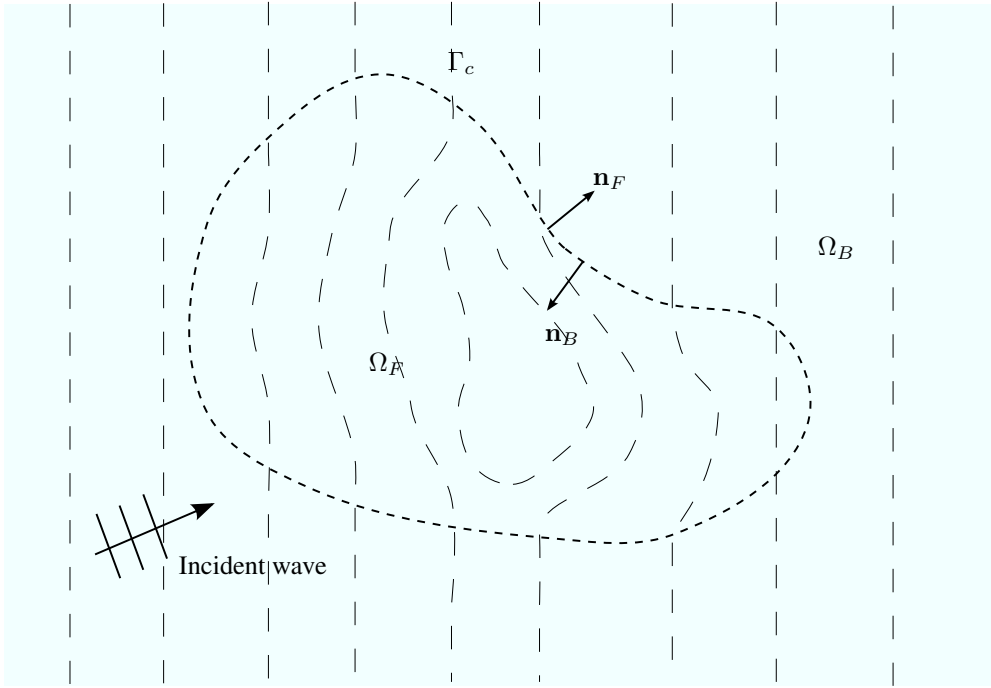
\caption{Illustration of the physical model with dashed lines representing iso-depth contours. A region with an arbitrary bathymetry is considered on $\Omega_{F}$ closed by $\Gamma_C$. On $\Omega_{B}$, the domain has a bathymetry with by straight and parallel contour lines} \label{fig:SchemeBEM_FEM}
\end{figure}
%
Formulations based on BEM and FEM offer different advantages that can be judiciously combined in order to palliate the drawbacks of both methods. The finite element approach of the MSE is able to model complex bathymetries while the boundary element formulation is able to satisfy automatically the Sommerfeld ration condition, being a perfect choice for modelling open and infinite domains.
A scheme of the coupling of the two methods proposed is showed in Figure \ref{fig:SchemeBEM_FEM}, with an incident wave train propagating over an infinite domain, that is divided into an inner and an outer regions. The inner region, $\Omega_F$, has an irregular bathymetry, while the outer region, named $\Omega_B$, presents parallel contour lines. 
The contour where both regions match is denoted $\Gamma_c$. A similar coupled BEM-FEM was proposed in \cite{Shaw1978}, but the bathymetry of the exterior domain was limited to be constant. In this work we include in the BSEM the fundamental solution derived by Cerrato et al. in \cite{Cerrato2016}, that makes it possible to consider a variable bathymetry in one direction, providing a perfect approach for the case under consideration.

The Spectral BEM-FEM coupling is carried out by a direct node to node matching condition along the boundary $\Gamma_c$. The vector containing the nodal values of the SEM domain, denoted by $\hat{\boldsymbol{\phi}}_{F} = [ \hat{\boldsymbol{\phi}}_{I} , \hat{\boldsymbol{\phi}}_{c}] $, that is divided into two parts, a subvector containing the interior nodes $\hat{\boldsymbol{\phi}}_{I}$ and another subvector containing the potentials located at the boundary, $\hat{\boldsymbol{\phi}}_{c}$. Then, combining the matrix forms obtained for the BSEM \eqref{eqn:BSEM_system} and SEM \eqref{eqn:SEM_system}, and considering the normal to $\Gamma_c$ positive towards the inner region, the non-homogeneous complete system can be written as
\begin{eqnarray}\label{BEM_FEM_coupling_sys}
\left[\begin{array}{ccc}
\hat{\mathbf{A}}_\mathrm{I} & \hat{\mathbf{A}}_\mathrm{c} & \hat{\mathbf{C}}\\
\mathbf{0}   & \hat{\mathbf{H}}   & -\hat{\mathbf{G}}
\end{array}\right] \left\{\begin{array}{c}
    \hat{\boldsymbol{\phi}}_\mathrm{I} 
    \\
    \hat{\boldsymbol{\phi}}_\mathrm{c} 
    \\
    \hat{\mathbf{q}}_\mathrm{c}
    \end{array}\right\} \!\!=  \left\{\begin{array}{c}
    \mathbf{0} 
    \\
    \hat{\boldsymbol{\phi}}_{in}
    \end{array}\right\},
\end{eqnarray}
where submatrices $\mathbf{A}_{\mathrm{I}}$ and $\mathbf{A}_{\mathrm{c}}$ collect the columns of $\mathbf{A}$ corresponding to $\hat{\boldsymbol{\phi}}_{I}$ and $\hat{\boldsymbol{\phi}}_{c}$ respectively. The incident wave is introduced into the formulation through the SBEM region considered as an scattering problem, where the SEM region is the scatter.

\section{Numerical results}\label{sec6}

The proposed formulation \eqref{BEM_FEM_coupling_sys} has been tested using different benchmark problems in order to demonstrate its high accuracy and convergence characteristics. The first example considers a wave train propagating over a constant water depth area using Neumann and Dirichlet boundary conditions to simulate the unbounded domain. This simple case allows us to evaluate the convergence and the behaviour of the solution of the BSEM and the SEM and compare them for the same problem. The second benchmark problem includes a variable bottom surface within the inner $\Omega_{F}$ region keeping a constant depth on the border$\Gamma_{C}$ with the outer region, testing the capabilities of the coupled spectral BSEM-SEM. The third example includes also a variable bathymetry in $\Omega_{B}$, being necessary to use a special fundamental solution for the BSEM.

\subsection{Convergence study of BSEM and SEM}\label{sec61}
\label{NE convergence} 
\begin{figure}[t]
\centering
\def\svgwidth{0.5\columnwidth}
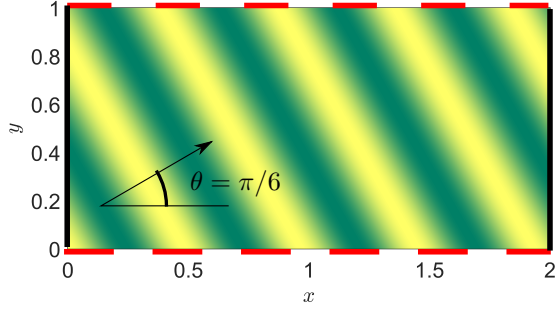
\caption{Wave train with an incident angle $\theta=\pi/6$ propagating through a rectangular domain with black solid lines representing Dirichlet boundary conditions and red dashed lines showing Neumann boundary conditions} 
\label{fig:Convergence Scheme}
\end{figure}

\begin{figure}[h]
\centering
\subfloat[]{\includegraphics[height=5.5cm,angle=0]{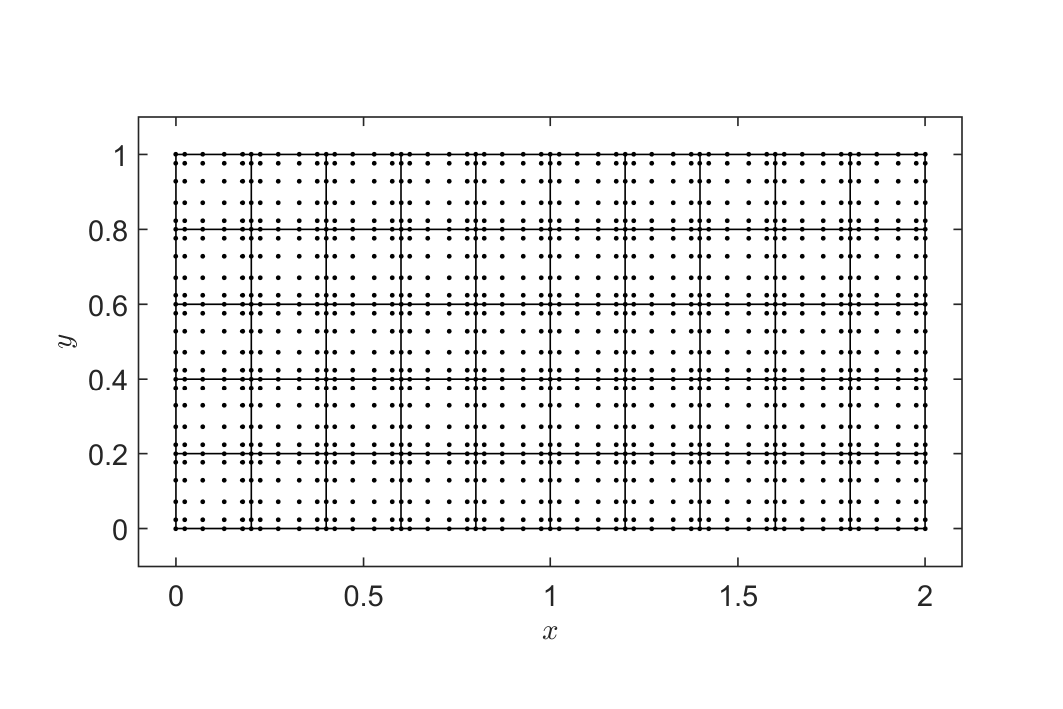}}
\subfloat[]{\includegraphics[height=5.5cm,angle=0]{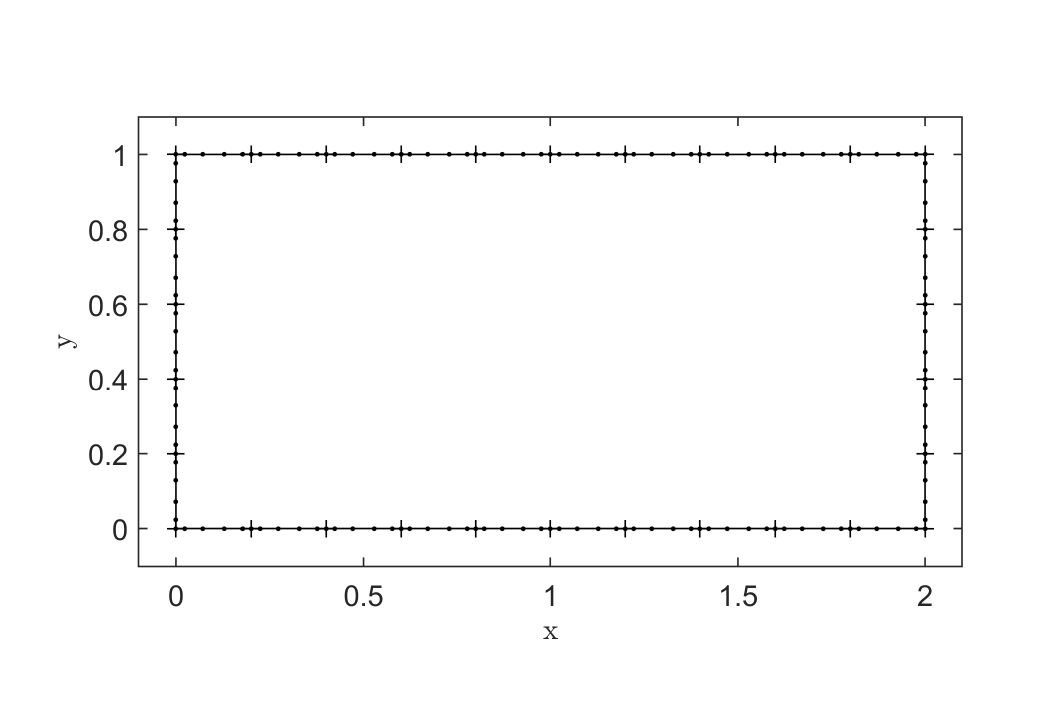}}
\caption{Spectral LGL finite (a) and boundary (b) element meshes for $h=0.2$ and $p=5$}
\label{fig:Convergence SEM_Mesh}
\end{figure}

\begin{figure}[t]
\centering
\includegraphics[width=0.6\columnwidth]{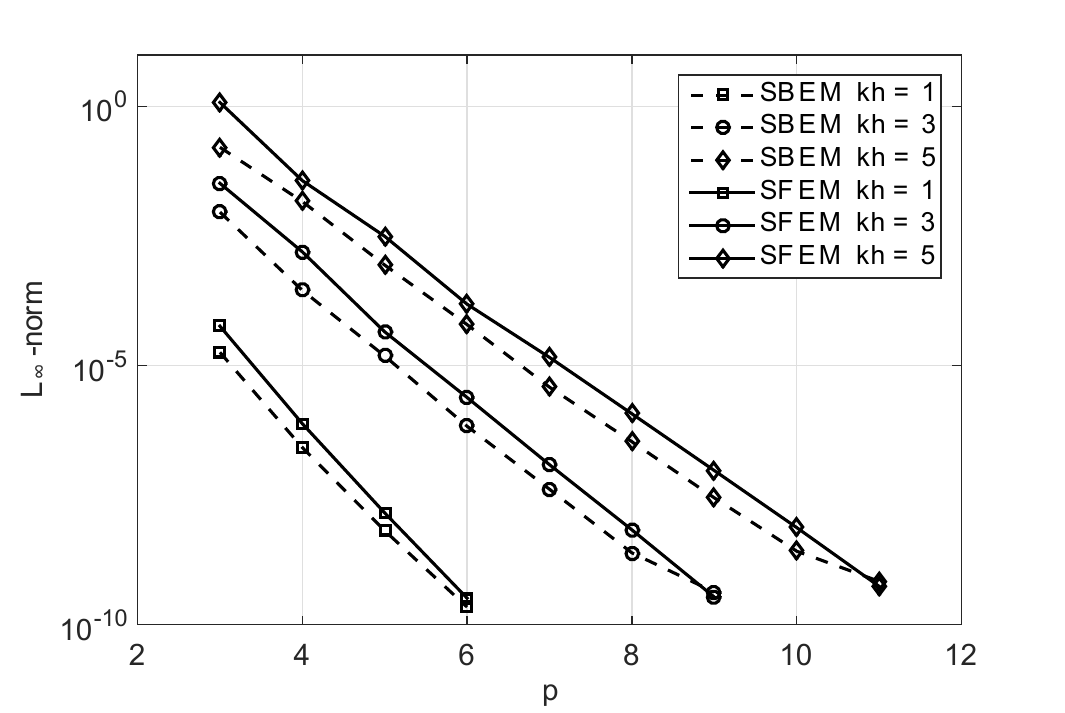}
\caption{A $p$-convergence study for different element size $h$} \label{fig:Convergence p-L2error}
\end{figure}
\begin{figure}[t]
\centering
\includegraphics[width=0.6\columnwidth]{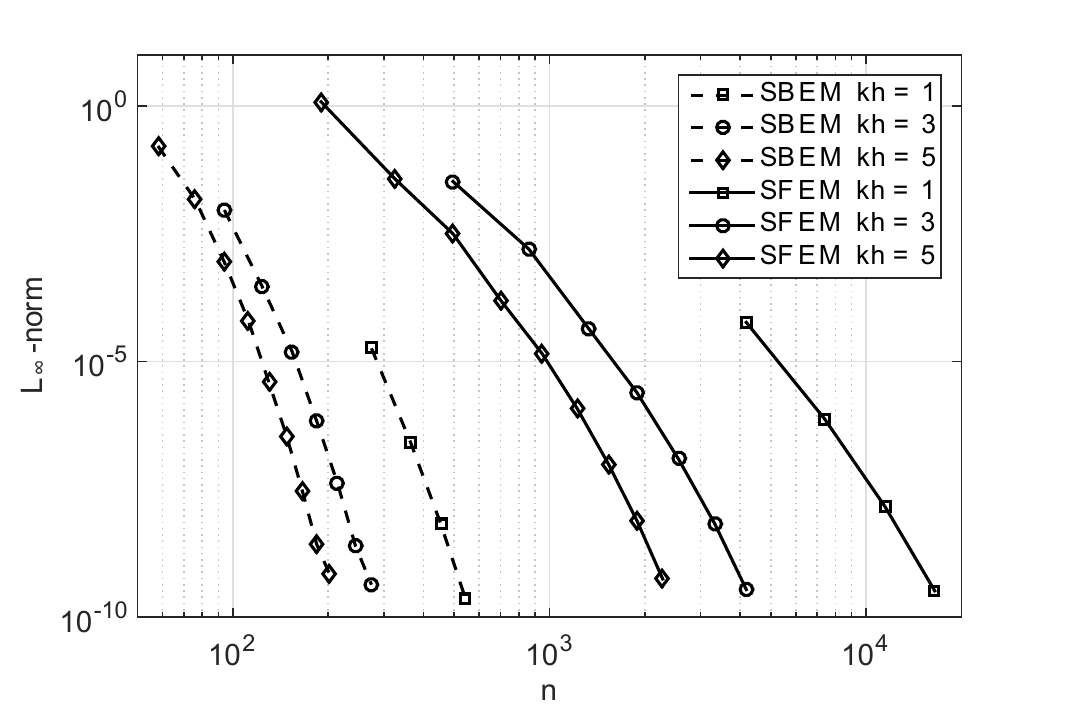}
\caption{Error evolution versus the number of degrees of freedom ($n$) for different number of elements per wave length} \label{fig:Convergence ndofs-L2error}
\end{figure}

This example considers a wave train propagating over an unbounded domain of constant-depth, as represented in Figure \ref{fig:Convergence Scheme}, being the solution $\hat{\phi}=e^{i\hat{k}(\cos\theta x+\sin\theta y)}$. The incident wave enters into the domain with an angle $\theta=\pi/6$. The wave number of the domain is $k=15$ and essential boundary conditions are applied on $x=0$ and $x=2$ (black solid lines in Fig. \ref{fig:Convergence Scheme}) while the known normal fluxes are applied on $y=0$ and $y=1$ (red dashed lines in Fig. \ref{fig:Convergence Scheme}).

The problem is solved using the SEM and the BSEM separately. Three different meshes are used with each method of elemet sizes $h=\left\lbrace 1/15,1/5,1/3\right\rbrace$ and therefore $kh=\left\lbrace 1,3,5 \right\rbrace$. To investigate the efficiency of the spectral LGL element formulations for BEM and FEM, p-convergence studies are done for every mesh, increasing the order $p$ of the nodal basis functions $L_m$. As example, the spectral finite element mesh for $h=0.2$ and $p=5$ is shown in Figure \ref{fig:Convergence SEM_Mesh} (a) and the corresponding spectral boundary element mesh is presented in Figure \ref{fig:Convergence SEM_Mesh} (b).

To analyze the convergence of the calculated potentials they are compared to the analytical solution. The $L_{\infty}$-error norm has been use for this purpose, which is defined as:
\begin{equation}
\label{eqn:NR Linfty definition}
||e||_{\infty}=\max\vert \boldsymbol{\hat{\phi}}_{\text{computed}} - \boldsymbol{\hat{\phi}}_{\text{exact}} \vert,
\end{equation}
where $\boldsymbol{\hat{\phi}}_{\text{computed}}$ is the vector that contains the computed solution at the nodes and $\boldsymbol{\hat{\phi}}_{\text{exact}}$ is the analytical solution. Analyzing Figure \ref{fig:Convergence p-L2error}, we can see that both methods (SEM and BSEM) present a standard spectral convergence, because the error norm decreases exponentially with the order of the approximation functions showing similar convergence slopes. Nevertheless, the BSEM presents a slighlty higher accuracy than the SEM. Another common feature of both methods is that the slope of the convergence curve increases when the element size $h$ decreases. Increasing the order of the nodal basis functions, is a better option to a certain limit, because it is known that the condition number of the system increases also linearly with the approximation order. 

The behaviour of the error-norm with the number of degrees of freedom ($n$) can be observed in Figure \ref{fig:Convergence ndofs-L2error}. As in Figure \ref{fig:Convergence p-L2error}, the number of elements remains constant, increasing only the order or the nodal basis function. Because only the boundary has to be discretized, the number of nodes needed to solve with the same accuracy the problem with the BSEM is much lower than those needed with SEM. However, regarding to the number of degrees of freedom, the main aspect that can be observed for both methods is again that, in every case, the strategy which seems more appropriate to achieve a fast convergence is to increase the element order instead of decreasing the element size. 

Observing the results of Figures \ref{fig:Convergence p-L2error} and \ref{fig:Convergence ndofs-L2error} we see that, a good choice is to use elements of size $h \approx 5/k$ with nodal basis functions of order $4\leq p \leq 8$ depending on the level of accuracy we need. 

Finally we can observe that both methods present similar accuracy when elements of the same order and size are used, so we can expect the same type of behaviour for a coupled finite and boundary SEM formulation.

\subsection{Wave scattering over a circular shoal}\label{sec62}
%
\begin{figure*}
\centering
\includegraphics[width=12cm]{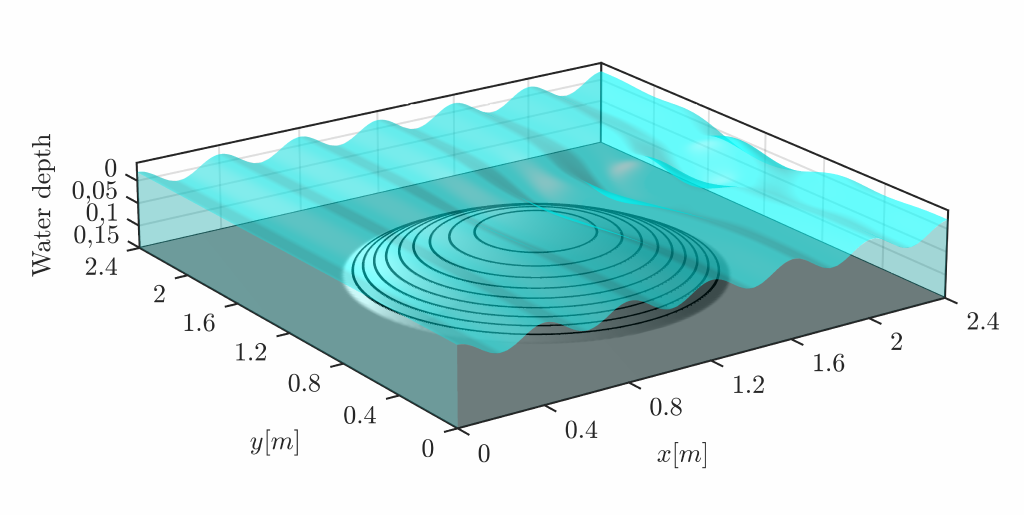}
\caption{
Water waves propagating over a circular shoal showing the effect downstream of an irregular seabed and its complete absorption by the coupled BSEM} \label{fig:Tanimoto scheme}
\end{figure*}
%


%
\begin{figure*}
\centering
\includegraphics[width=\textwidth]{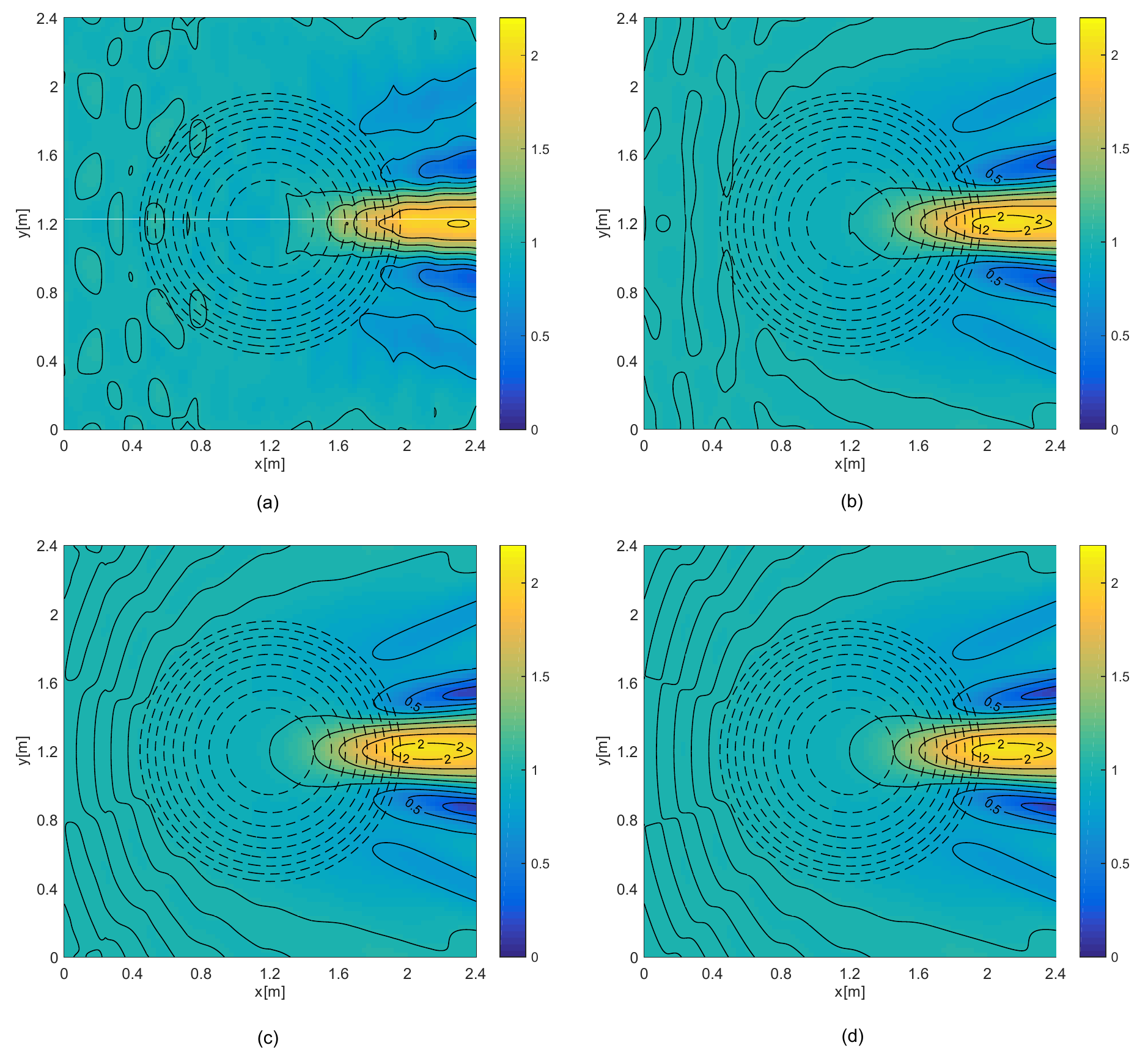}
\caption{
Water-wave propagation over a circular shoal. Absolute value of the dimensionless wave height field over the circular shoal (i.e. dashed lines) for different element polynomial orders: (a) $p=3$, (b) $p=4$, (c) $p=5$ and (d) $p=15$} \label{fig:Tanimoto_errors_surfs}
\end{figure*}
%


\begin{figure}[!t]
\centering
\subfloat[]{\includegraphics[height=5.05cm,angle=0]{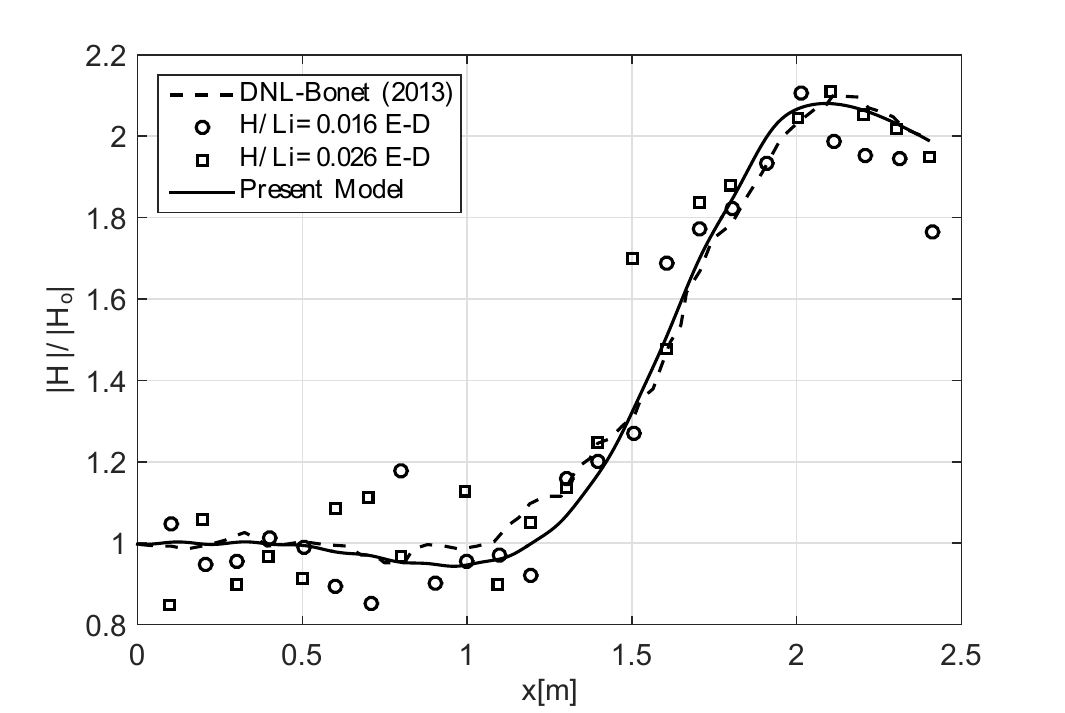}}
\subfloat[]{\includegraphics[height=5.05cm,angle=0]{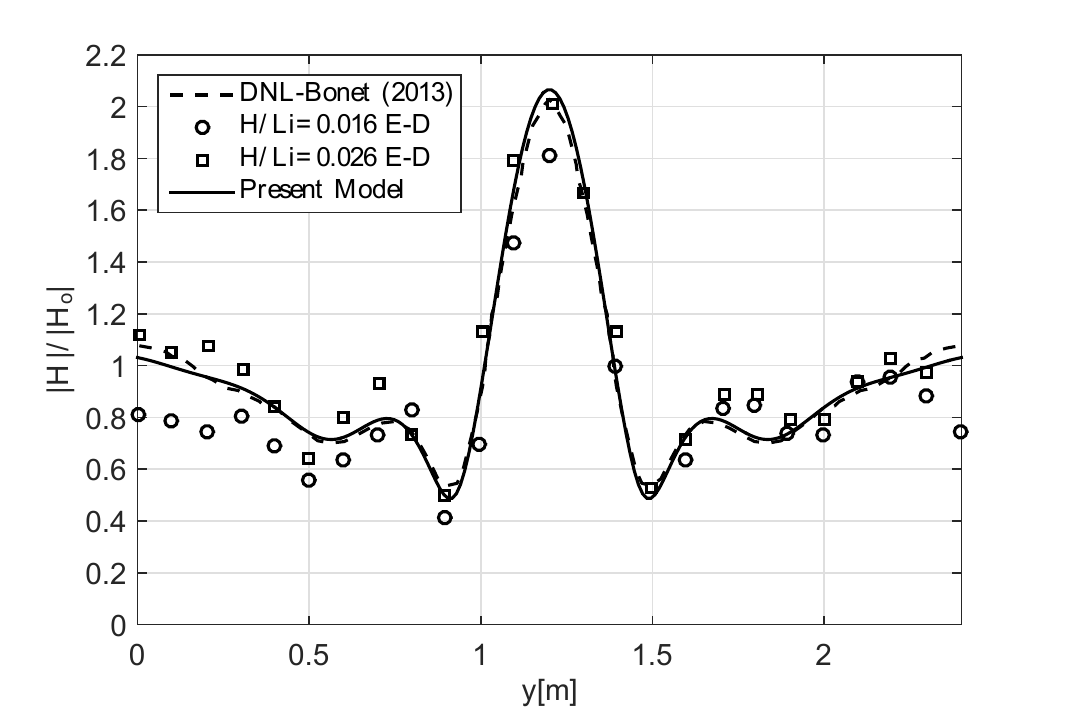}}\\
\subfloat[]{\includegraphics[height=5.05cm,angle=0]{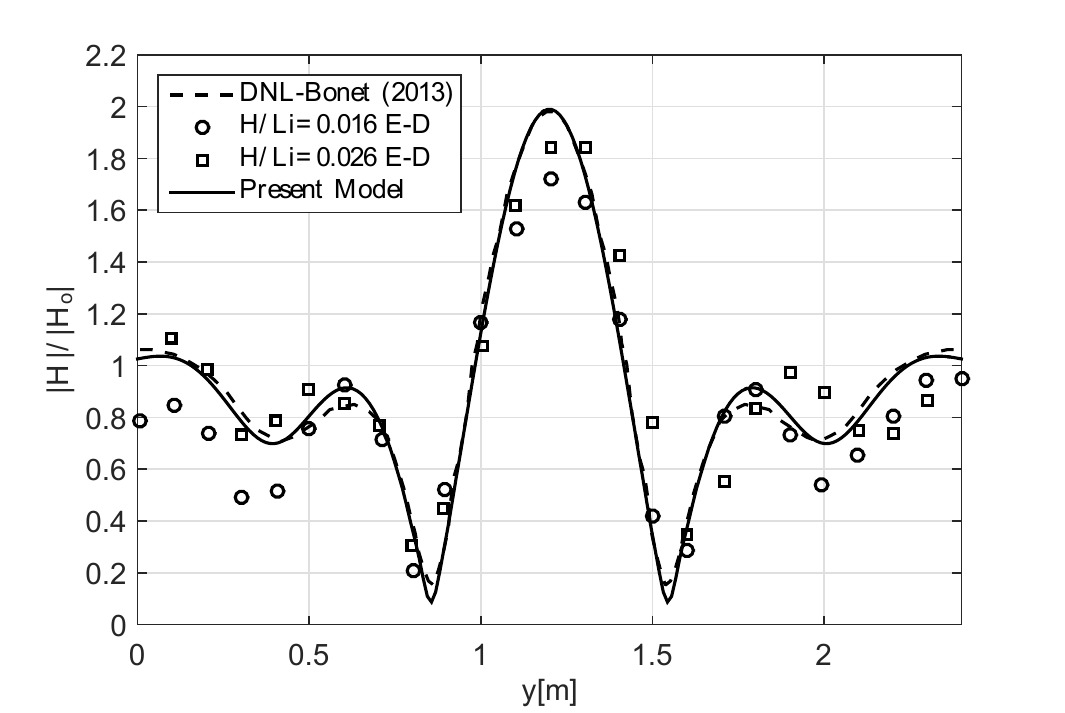}}
\caption{Normalized wave height over a circular shoal. Comparison of the present coupled spectral BEM-FEM solution (solid line) with experimental data from \cite{Ito1972} and the solution presented in \cite{Bonet2013} for three different sections: (a) section $y=1.2$, (b)  section $x=2$ and (c)  section $x=2.4$}
\label{fig:Tanimoto sections}
\end{figure}

\begin{figure}[!t]
\centering
\subfloat[]{\includegraphics[height=4.5cm,angle=0]{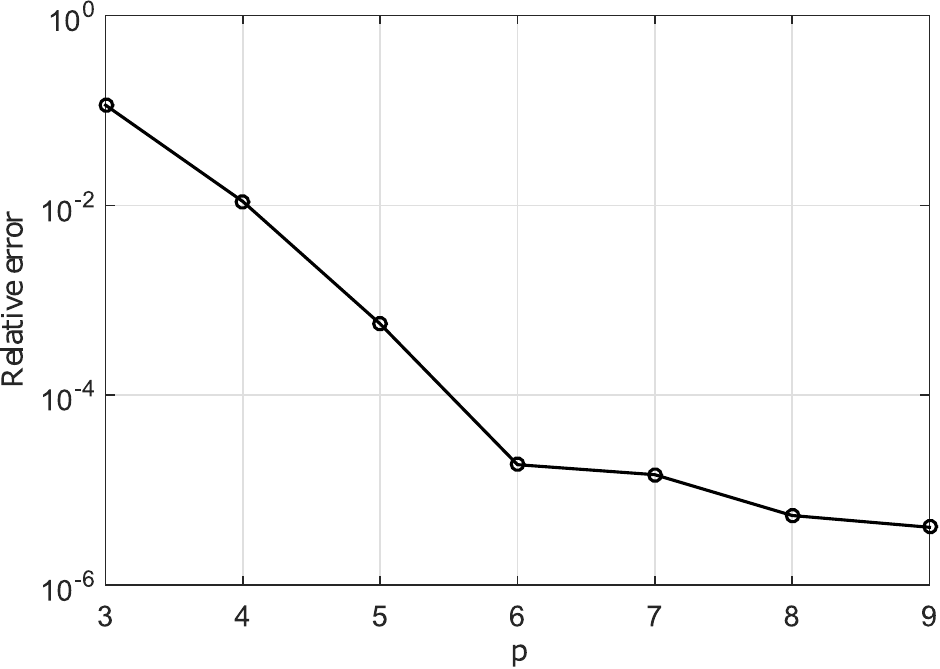}}
\hspace{0.5cm}
\subfloat[]{\includegraphics[height=4.5cm,angle=0]{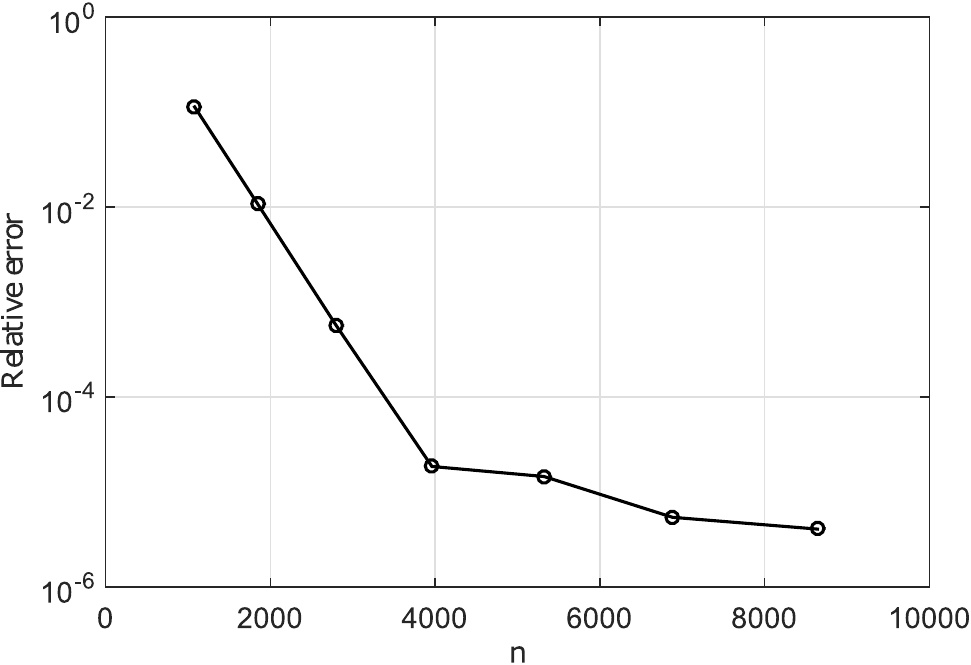}}
\caption{Water-wave propagation over a circular shoal. A $p$-convergence study for a mesh with $10\times10$ spectral finite elements and $40$ spectral boundary elements coupled on the boundary. (a) Relative error versus polynomial order ($p$). (b) Error evolution with the number degrees of freedom ($n$)} \label{fig:Tanimoto_errors}
\end{figure}

Like in the previous problem, a wave train propagating in an infinite domain is considered but, this time, we include on the seabed a parabolic shoal. The water depth is defined by:
\begin{equation}
h_w(r)=\left\lbrace
\begin{array}{l l}
0.1\left(\dfrac{r}{0.8}\right)^2+0.05, & r<0.8,
\\
0.15, & r\geq 0.8.
\end{array}
\right.
\end{equation}
where $r$ is the distance to the center of the shoal: $r=\sqrt{(x-x_c)^2+(y-y_c)^2}$, located at $(x_c=1.2,y_c=1.2)$, and the water depth takes values between $h_w=0.15$ and $h_w=0.05$. The internal domain is a square located in $(0\leq x \leq 2.4, \: 0\leq y \leq 2.4$, with an infinite region of constant depth around it. 

A wave train of period $T=0.511s$ and different wave heights have been simulated. A scheme of the model, where the wave height has been exaggerated, is showed in Figure \ref{fig:Tanimoto scheme}. The relation between the water depth and the wave period indicates that waves are propagating over intermediate waters, where the MSE is considered to be applicable. 

The partition of the model into an internal region enclosing the irregular bathymetry and another external region with regular bathymetry, makes of this case an excellent benchmark to test the capabilities of the coupled spectral BSEM-SEM model. As it was described in previous sections, SEM allows to model regions including arbitrary bathymetries, but, if we use only this method, special techniques such as PML are required to model correctly open boundaries avoiding spurious reflection. The BSEM is used instead to model the infinite domain, fulfilling automatically the Sommerfeld radiation condition. Considering that the interface $\Gamma_c$ between both regions is located over constant depth, the fundamental solution used for the BSEM is the standard Green's function for the two dimensional Helmholtz problem (more details are given in the appendix \ref{apendiceA}). This is an important element because it allows us to study the problem for different mesh configurations without introducing additional sources of error related to the absorption of the outgoing waves.

To approximate the wave field, the domain represented in Figure \ref{fig:Tanimoto scheme} is discretized using a rectangular mesh of $10\times10$ spectral finite elements of the same length. The SEM mesh is surrounded by a perfectly matching BSEM mesh composed by 40 elements. Then, the problem is studied as an scattering problem, with an incident wave field going in the $x$-direction. The selection of the element size was realized following the conclusions derived from the convergence analysis. From the dispersion relation we know that the wave number takes values between $15.69$ $rad/m$ and $20.15$ $rad/m$, with values of the dimensionless parameter $kh$ between 3.77, at the center of the shoal, and 4.84 over the flat region. 

To analyse the behaviour of the coupled formulation, the problem is solved using different order of the basis approximation functions and compare the results using the following definition for the global error norm: 
\begin{equation}
\label{eqn:Ito-Relative error}
\text{Relative error} = \left\vert \dfrac{\int_{\Omega} (\hat{\phi}-\hat{\phi}_{R}) \, d\Omega}{\int_{\Omega} \hat{\phi}_{R} \, d\Omega} \right\vert,
\end{equation}
where the computed solution $\hat{\phi}$ is compared with a reference solution $\hat{\phi}_{R}$.

Unfortunately there is no analytical solution for this problem, but experimental data was obtained in 1972 by Ito and Tanimoto \cite{Ito1972}, that designed this experiment to test refraction and diffraction numerical models. In order to select a reference solution $\hat{\phi}_R$, the problem has been solved repeatedly increasing the element order until convergence with $p=15$.
The convergence of the results is shown in Figure \ref{fig:Tanimoto_errors_surfs}, representing the normalized wave-height ($H/H_o$) for different element orders $p=\left\lbrace3,4,5,15\right\rbrace$. The solution varies drastically when the order of element ischnaged from $3$ to $4$, and small improvements are observed changing the order from $4$ to $5$, indicating that the reference solution has been practically reached. From order $5$ to $15$ the solution does not change significantly, but the accuracy of the solution continues increasing. Then, the normalized wave field calculated with $p=15$ has been chosen as the correct solution. 
To validate our results, some sections of the wave field are compared with the experimental data presented in \cite{Ito1972} and the numerical results from Bonet \cite{Bonet2013}, where standard finite elements are used in combination with DNL boundary conditions. 
This comparison is shown in Figure \ref{fig:Tanimoto sections}. Section (a) represents the profile of the normalized wave height along the $x$-axis at $y=1.4$, section (b) gives the profile along the $y$-axis at $x=2$ and finally section (c) is located at $y=2.4$. The experimental data represented in the figure, correspond to two different wave heights: $H=0.0064m$ and $H=0.0104m$. All the profiles show a good agreement with the solutions from other authors. We remark the smoothness of our solution, which is not typically found on other methods based on a low order approximation of the solution due to the pollution effect. The use of high-order polynomials and the derivability of the solution works really well diminishing the intrinsic error of discrete methods. 

Finally, the error norm is used to analyse de behaviour of the solution. Figure \ref{fig:Tanimoto_errors} represents the relative error versus the element order and the number of degrees of freedom. The convergence is similar to the previous example, but this time is truncated at $10^{-4}$ and $10^{-5}$. The origin of this blocking in the convergence is in the difficulty to integrate the $\hat{k}(x,y)$-function, with a presenting non-continuous first derivative on the edge of the circular shoal. However, even considering an inhomogeneous medium, the spectral convergence of the method is still preserved, providing high accuracy and smoothness of the solution.

\subsection{Wave scattering over a sloping bottom with an elliptical shoal}\label{sec63}

%
\begin{figure*}[t]
\centering
\includegraphics[width=\textwidth]{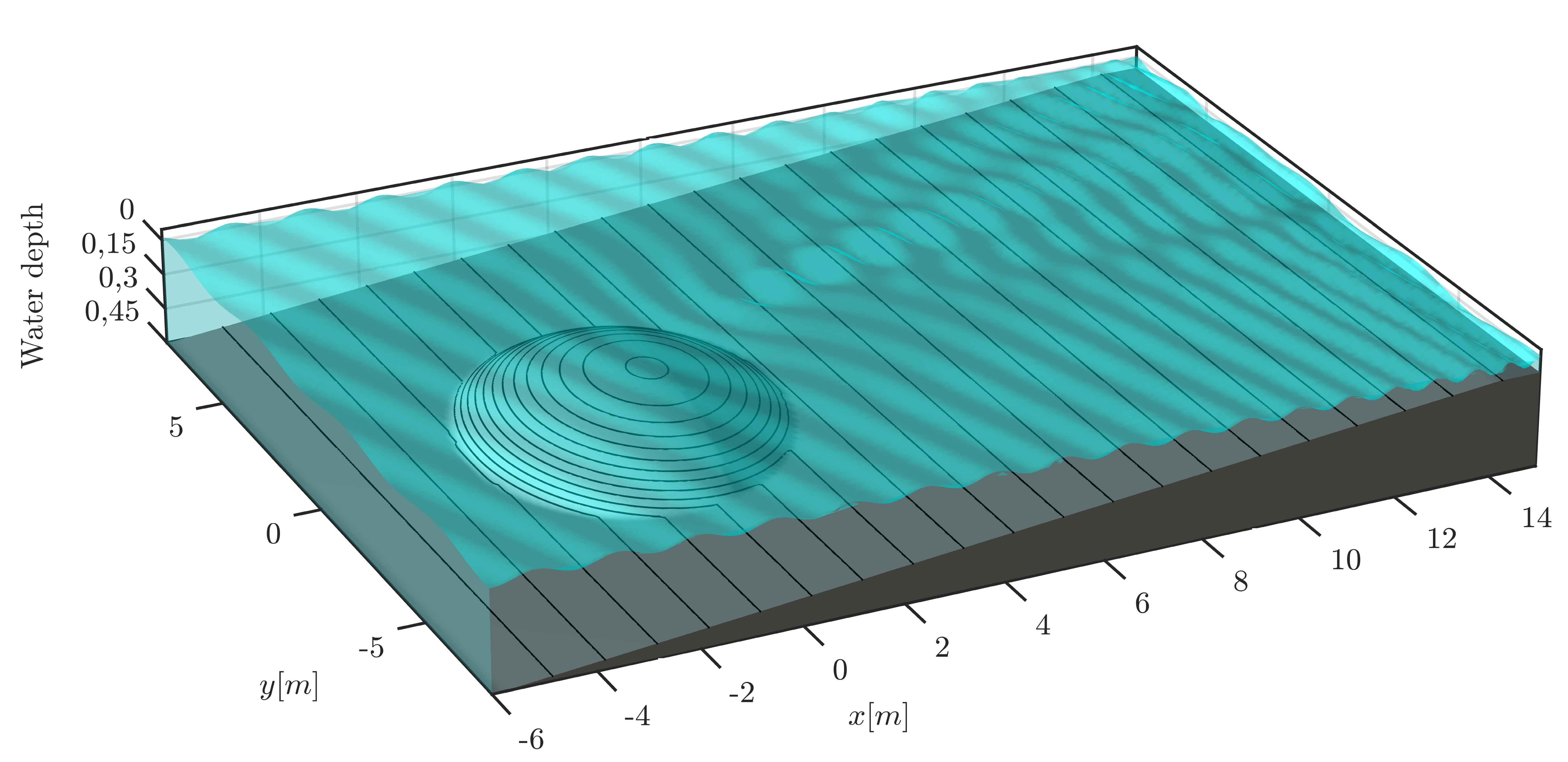}
\caption{Water-wave propagation problem over a sloping bottom with an elliptical shoal} \label{fig:Elip shoal scheme}
\end{figure*}
\begin{figure*}[ht!]
\centering
\subfloat[$p=3$]{\includegraphics[width=0.6\textwidth,angle=0]{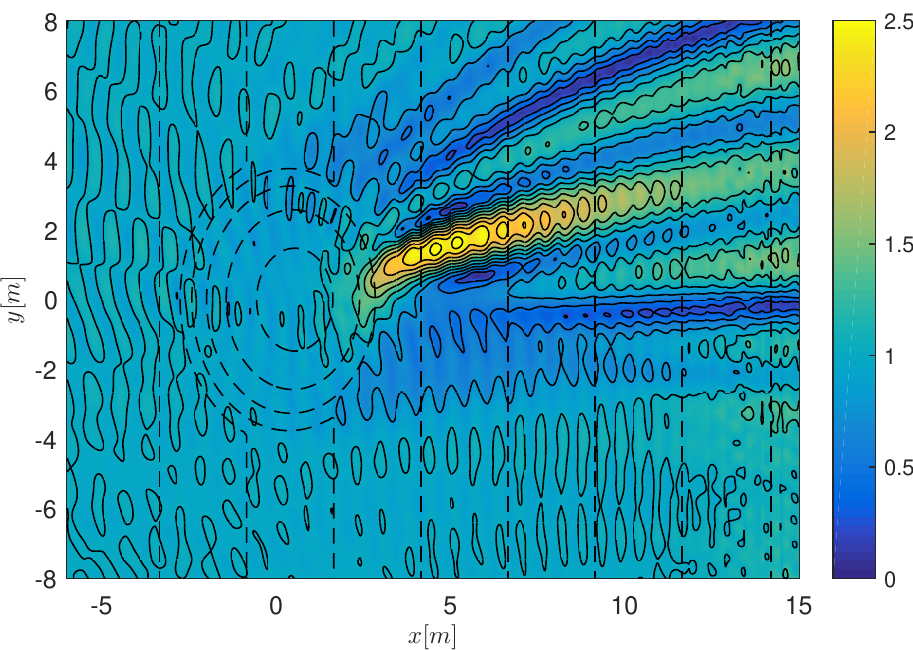}}
\vspace{0.1cm}
\subfloat[$p=4$]{\includegraphics[width=0.6\textwidth,angle=0]{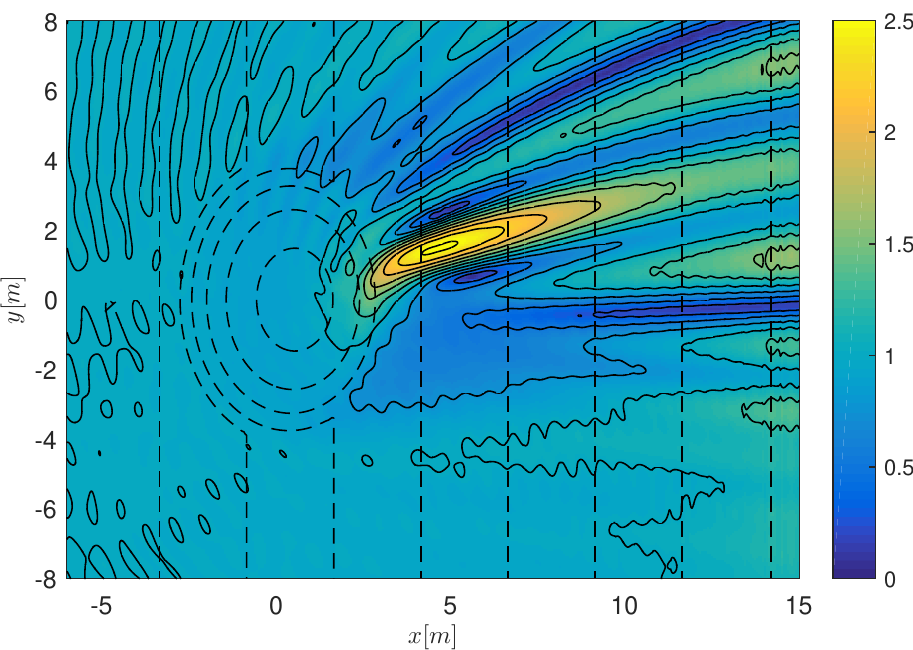}}
\vspace{0.1cm}
\subfloat[$p=6$]{\includegraphics[width=0.6\textwidth,angle=0]{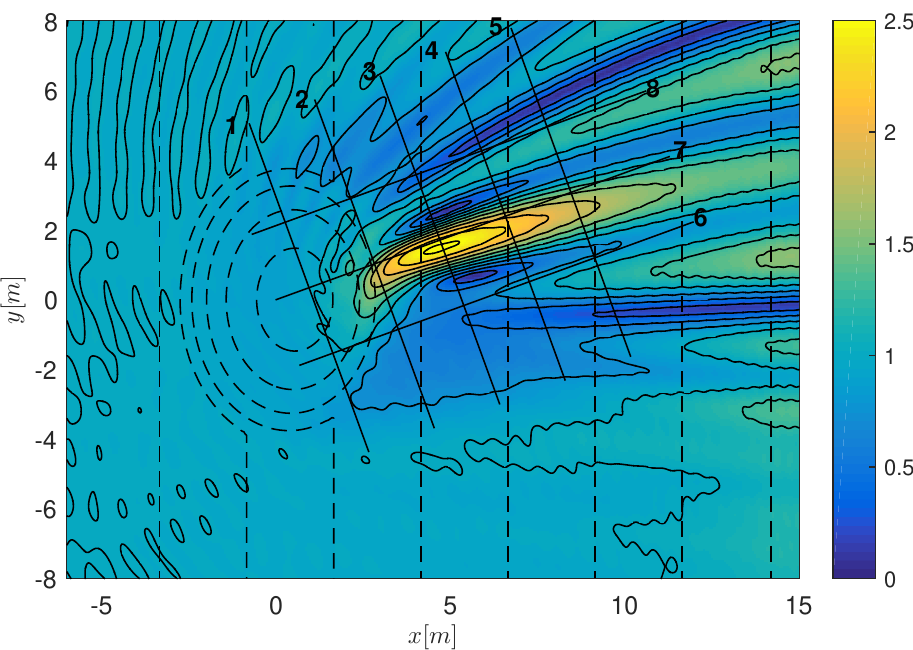}}
\caption{Elliptical shoal problem. Absolute value of the dimensionless wave height field over the elliptical shoal for: (a) $p=3$, (b) $p=4$ and (c) $p=6$} \label{fig:Elip shoal solutions}
\end{figure*}

The elliptic shoal over a constant slope bed is a common numerical experiment which has been widely used to demonstrate the behabiour and capabilities of different wave theories and numerical techniques for water wave \cite{Li1994,Berkhoff1982,Belibassakis2001}. A geometrical description of the problem is depicted in Figure \ref{fig:Elip shoal scheme}, where depth contours and water surface are represented. The seabed is dominated by a constant slope of $2\%$ starting from $0.45m$-depth up to $0.05m$-depth. Over this slope, an elliptic shoal is superimposed. The depth function is defined in metres as $h_w(x,y)=h_i+h_s$ where:
\begin{equation}
h_i(x,y)=\left\lbrace
\begin{array}{l l}
0.45,& x<-5.85, \\
0.45-0.02(5.85+x), & -5.85 \leq x \leq 14.15, \\
0.05, & x>14.15,
\end{array}
\right.
\end{equation}
is the expression for the constant slope component and the superimposed shoal produces a disturbance height 
\begin{equation}
h_s(x,y)=0.3-0.5\sqrt{1-\left(\dfrac{x}{3.75}\right)^{2}-\left(\dfrac{y}{5}\right)^{2}} \hspace{0.2cm} \text{in} \hspace{0.2cm} \Omega_{s},
\end{equation}
where $\Omega_{s} = \{(x,y) : \:(x/3)^{2}+(y/4)^{2} \le 1$.
Experimental data for this problem are provided in \cite{Berkhoff1982}, with an incident wave train of period $T=1s$ and height of $0.01058m$, the direction of the waves is $\theta=20^o$ relative to the $x-axis$. The maximum slope of the seabed is on the front of the elliptical shoal, where it raises up to $18\%$ and the shallowness ratio ($h_w/\lambda$) varies from $0.3$, on the deepest zone, to $0.074$, on the shallow part, meaning that waves are propagating over intermediate waves. These geometric parameters of the bottom profile indicate that the MSE is adecuate to model the problem, nevertheless the amplitude of the incident wave used to obtain the experimental data on \cite{Berkhoff1982} falls outside the limit of linear theory on the shallow zone and if one wants to obtain a better approach of the physical reality it should be considered to use non-lineal models. 

The aim of solving this case by the proposed coupled BSEM-SEM formulation, is to analyse its behaviour when the bathymetry varies also in the external infinite region. To make it possible, a specialized fundamental solution for variable water depth, like the one presented in \cite{Cerrato2016}, has to be included.

Based on the experience from convergence analyses, a regular mesh of spectral finite elements is used to discretize the $xy$-domain, with 40 divisions along the $x$-axis and 30 divisions along the $y$-axis, which makes a total of 1200 spectral finite elements. The external interface is meshed using 140 spectral boundary elements perfectly matching the spectral finite element discretization on the boundary. In this regular mesh the $kh$ coefficient remains between 2.2 and 4.7.

The absolute value of the normalized wave-heigh has been obtained for different element orders, starting at $p=3$. Some of the calculated solutions are shown in Figure \ref{fig:Elip shoal solutions}, with $p=3$ (a), $p=4$ (b) and $p=6$ (c). The results reveal a fast convergence, being the difference between orders 3 and 4 remarkable, while very small changes are observed on the solution from order 4 to 6. Order 6 is enough to represent the correct solution with a crest behind the shoal with two important amphidromical points, corresponding to low amplitude areas, as described in \cite{Berkhoff1982,Belibassakis2001}.

\begin{figure*}[ht]
\centering
\includegraphics[width=\textwidth]{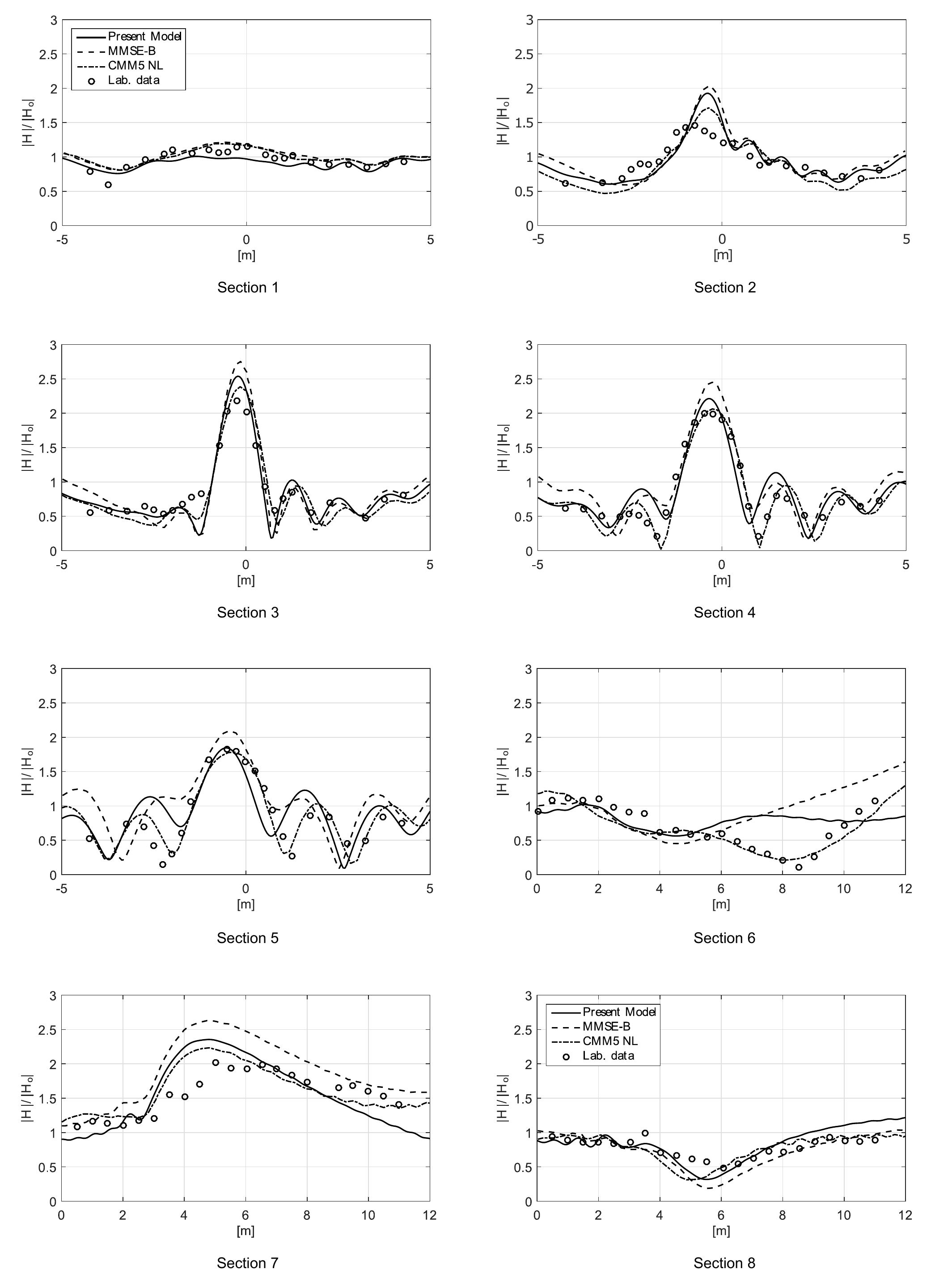}
\caption{Elliptical shoal problem. Comparison of the coupled spectral BEM-FEM solutions (solid line) with the MMSE (dashed line), a Coupled-Mode Model (dashed-point line) \cite{Belibassakis2001} and the experimental data from \cite{Berkhoff1982} for different sections} \label{fig:Elip shoal sections}
\end{figure*}

Figure \ref{fig:Elip shoal sections} shows the sections which are marked on Figure \ref{fig:Elip shoal solutions} (c). In these profiles, the normalized wave-height field calculated using the proposed coupled BSEM-SEM formulation (solid line) with $p=6$ is compared to the MMSE (dashed line), the values from a Coupled-Mode Model (dashed-point line) \cite{Belibassakis2001} and the experimental data from \cite{Berkhoff1982}. A good agreement with other solutions can be observed. Nevertheless, as expected, the mild-slope models cannot approximate experimental data as well as the coupled-mode model does, because the MSE does not exactly preserve mass conservation due to the bottom curvature and also because slope-squared terms have been neglected.

\section{Summary and conclusions}\label{sec7}

Realistic wave propagation problems modelled by the MSE need to include the effect variable bathymetries, the modelization of big areas and the correct approach of open and infinite regions. At the same time, standard FEM formulations based on low order polynomial approximations require a big number of degrees of freedoms, with an increasing number nodes per wave length as the wave number increases due to the pollution effect. Considering all these aspects and in order to provide a better approach, a coupled boundary spectral element method (BSEM) and finite spectral element method (SEM) formulation has been proposed. This formulation is based on a pseudospectral approximation, considering elements with a Legendre-Gauss-Lobatto (LGL) distribution of nodes. 

A complete formulation of the BSEM has been derived for the MSE. This new method take advantages of the specific virtues of the spectral approach with LGL elements and the Boundary Integral Formulation (BIE) of the Helmholtz problem with variable wave number. One of the difficulties of BIE formulations is in the integration of singular kernels, which increases the computational cost with the order of the approximation polynomials. Therefore, a regularization technique has been proposed, with an analytic integration of the singularity and a regular integral of a smooth function that can be computed by using a standard Gaussian quadrature.

The results of this work show that a pseudospectral approach of the solution provides a similar convergence in BSEM and SEM, achieving very accurate solutions, reducing the pollution effect drastically with relatively few degrees of freedom, compared to classical discrete approaches of the problem. 

The coupled BSEM-SEM strategy can cover a wide range of different bathymetries. The SEM is used to model internal regions, usually the zone of interest. SEM combines all the positive attributes of the classical FEM approach toguether with the accuracy of spectral methods, being able to reproduce water wave propagation over variable bathymetries with high accuracy. In order to reproduce natural boundary conditions, the internal region is connected to an outer infinite region modelled by the BSEM, which is based on the Boundary Integral Equation, automatically fulfiling the Sommerfeld radiation condition. For the BSEM implementation, the complete fundamental solution presented in \cite{Cerrato2016} for variable bathymetries has been also included, increasing the capabilities of the coupled BSEM-SEM modelization. This fundamental solution allows to reproduce more realistic bathymetries, by locating the coupling frame in regions where the only restriction is that the bathymetry lines have to be straight and parallel contour lines.

Two classical water-wave benchmark problems has been solved using the coupled formulation for the MSE, providing an accurate solution of the wave field and the water profiles. 

As a future work we consider the possibility of extending this method to other water wave formulations, such as the MMSE, or to include the effect of wave breaking or currents. Also, the interaction of waves with structures \cite{Luis2015} is an interesting topic for wave energy generation.

\section{Acknowledgments}
This work was supported by the \emph{Ministerio de Ciencia e Innovaci\'{o}n} (Spain) through the research project DPI$2010$-$19331$.

\clearpage


\section{Appendix: Fundamental solutions}\label{apendiceA}
The fundamental solution for the 2D-Helmholtz problem with variable wave number $\hat{k}(x)$ in one direction, is obtained from the solution of:
\begin{equation}
\label{eqn:Green}
\nabla^{2}\psi+\hat{k}(x)^{2}\psi + \delta(\xbold-\xbold') = 0 
\end{equation}
with $\xbold'$ the collocation point and considering the Sommerfeld ration boundary condition at infinity. For the special case where $\hat{k}$ is constant the fundamental solution is given by the following expressions:
\begin{equation}
\psi(\xbold,\xbold',\hat{k})=  \dfrac{i}{4} H_{0}^{(1)}(\hat{k}r) , \quad
\psi_{,r}(\xbold,\xbold',\hat{k}) = - \dfrac{i}{4} \hat{k} H_{1}^{(1)}(\hat{k}r),
\end{equation}
where $H_{0}^{(1)}$ and $H_{1}^{(1)}$ are Hankel functions of the first kind of order zero and one, respectively, and $r$ is the distance to the collocation point.

For variable bathymetries changing only in one direction ($x$) and with a wave number function described by a continuous function of the form:
\begin{equation}
\hat{k}(x)=\left\lbrace
\begin{array}{l l}
\hat{k}_a,& x<a, \\
\hat{k}_b(x), & a \leq x \leq c, \\
\hat{k}_c, & x>c,
\end{array}
\right.
\end{equation}
it is possible to apply a Fourier transform ($\mathcal{F}(\bullet)$) in the $y$-direction to \eqref{eqn:Green}, leading to the following one-dimensional problem:
\begin{align}
\label{eqn:1Dvarphi in 2}
& \varPsi_{,xx} + \kappa^{2}(x) \varPsi + \delta(x - x') = 0 \quad \text{in } x \in [a,c] \\
\label{eqn:bc-varphi-x}
& \left\{\begin{array}{lll}
\varPsi_{,x} + i \alpha(\xi) \varPsi  =  0 & \text{in} & x=a \\
\varPsi_{,x} - i \beta(\xi) \varPsi  =  0 & \text{in} & x=c \\
\end{array}\right.
\end{align}
where $\varPsi(x,x';\xi) = \mathcal{F}(\psi(\xbold,\xbold'; \hat{k}))$ and $\kappa^{2}(x)=\hat{k}_b^{2}(x)-\xi^{2}$ are the transformed velocity potential and wave number. The coefficients $\alpha$ and $\beta$ are: $\alpha(\xi)=(\hat{k}_a^2-\xi^2)^\frac{1}{2}$ and $\beta(\xi)=(\hat{k}_c^2-\xi^2)^\frac{1}{2}$. This one-dimensional problem can be solved numerically using different techniques, like for example the SEM described in Section \ref{sec4}. 

Following \cite{Cerrato2016}, the complete kernel of the boundary element formulation for linear water waves propagating over a variable bathymetry can be calculated by solving the following integrals coming from the application of the inverse Fourier transform ($\psi=\mathcal{F}^{-1}(\varPsi)$) to the solution of (\ref{eqn:1Dvarphi in 2},\ref{eqn:bc-varphi-x}), obtaining:
%
%
\begin{equation}
\label{eqn:IFT-phi}
\psi(\xbold,\xbold'; \hat{k}) = \dfrac{1}{\pi} \int_{0}^{\infty} \varPsi (x,x';\xi) \cos{( y \xi)} \: d\xi
\end{equation}
for the velocity potential and
\begin{align}
\label{eqn:x-derivative-psi}
\psi_{,x}(\xbold,\xbold'; \hat{k}) &=
\dfrac{1}{\pi} \int_{0}^{\infty} \varPsi_{,x} (x,x';\xi) \cos{( y \xi)} \: d\xi, 
\\
\label{eqn:y-derivative-psi}
\psi_{,y}(\xbold,\xbold'; \hat{k}) &=
\dfrac{1}{\pi} \int_{0}^{\infty} - \xi  \varPsi (x,x';\xi) \sin{( y \xi)} \: d\xi
\end{align}
for the derivatives of the Green's function. Details about the numerical computation of these integrals can be found in \cite{Cerrato2016}.

\bibliographystyle{unsrt}  



\end{document}